\documentclass[a4paper,12pt,twoside]{amsart}
\usepackage{amsmath,amsthm,amsfonts}
\usepackage[T1]{fontenc}
\newtheorem{thm}{Theorem}[section]

\newtheorem{lem}[thm]{Lemma}

\newtheorem{prop}[thm]{Proposition}
\newtheorem{proposition}[thm]{Proposition}

\theoremstyle{definition}

\newcommand{\R}{{\mathbb{R}}}

\newcommand{\E}{{\mathbb{E}}}

\newcommand{\F}{{\mathbb{F}}}
\newcommand{\Z}{{\mathbb{Z}}}
\newcommand{\N}{{\mathbb{N}}}
\newcommand{\T}{{\mathbb{T}}}

\newcommand{\Cal}{\mathcal}
\newcommand{\cal}{\mathcal}
\newcommand{\cH}{\mathcal{H}}

\def \Ker {{\rm Ker \,} }
\def \mod {{\rm \ mod \,} }
\def \diam {{\rm diam} \,}

\def\g{{ \mathbf g}}
\def\j{{\underline j}} \def\k{{\underline k}}
 \def\el{{\underline \ell}}
\def\t{{\underline t}} 
\def\0{{\underline 0}}  \def\1{{\underline 1}}
 \def\g {{\underline g}}
\def\a{{\underline a}} \def\b{{\underline b}}
\def\c{{\underline c}}
 \def\p {{\underline p}}
 
\def\v{{\underline v}} 
\def\s{{\underline s}} \def\x{{\underline x}}

\def\uA{{\underline A}}
\def\uR{{\underline R}}

\def\al {{\underline \alpha}} \def\be {{\underline \beta}}

\def \and{\text{ and }}

\def \stm0{{\setminus \{\0\}}}

\def \eop{\qed}

\def \proof {\vskip -3mm {{\it Proof}. }}

\def \Log {{\rm Log}}

\begin{document}
\title[$\Z^d$-actions on subgroups of $\F_p^{\Z^d}$]
{Almost mixing of all orders and CLT \\
for some $\Z^d$-actions on subgroups of $\F_p^{\Z^d}$}

\date{18 September 2016}
\author{Guy Cohen and Jean-Pierre Conze}
\address{Guy Cohen, \hfill \break Dept. of Electrical Engineering,
\hfill \break Ben-Gurion University, Israel}
\email{guycohen@bgu.ac.il}
\address{Jean-Pierre Conze,
\hfill \break IRMAR, CNRS UMR 6625, \hfill \break University of
Rennes I, Campus de Beaulieu, 35042 Rennes Cedex, France}
\email{conze@univ-rennes1.fr}


\begin{abstract}
For $\mathbb{N}^d$-actions by algebraic endomorphisms on compact abelian groups,
the existence of non-mixing configurations is related to "$S$-unit type" equations and plays a role in limit theorems for such actions.

We consider a family of endomorphisms on shift-invariant subgroups of $\mathbb{F}_p^{\mathbb{Z}^d}$ and show that there are few solutions
of the corresponding equations. This implies the validity of the Central Limit Theorem for different methods of summation.
\end{abstract}

\maketitle
\tableofcontents

\section*{\bf Introduction}

Let $G$ be a compact abelian group endowed with its Haar measure $\mu$. If $T_1^{\ell_1}, ..., T_d^{\ell_d}$, $d \geq 1$,
are commuting algebraic automorphisms or surjective endomorphisms of $G$, they generate a $\Z^d$ or $\N^d$-action on $G$:
$\el \to T^\el =T_1^{\ell_1} ... T_d^{\ell_d}$.

Given an "observable" $f: G \to \R$ with some regularity, one can investigate
the statistical behavior of the random field $(T^\el f)_{\el \in \N^d}$, in particular the following limits (in distribution with respect to $\mu$):
\hfill \break- ergodic sums for a sequence $(D_n)$ of sets in $\N^d$:
$\lim_n \, |D_n|^{-\frac12} \  \sum_{\el \in D_n} T^\el f,$
\hfill \break- ergodic sums along a random walk $Z_n = Y_0 + ... + Y_{n-1}$ on $\Z^d$ or $\N^d$:
$$\lim_n \, a_n^{-1} \, \sum_{0 \leq k < n} T^{Z_k(\omega)} \, f,  \text{ for a.e. fixed } \omega,
\text{ where } (a_n) \text{ is a normalizing sequence}.$$
The connected case was considered in \cite{CohCo15}. Here we are interested in non connected groups $G$. More precisely we consider
in Section \ref{shiftInv} some commutative actions by endomorphisms or automorphisms on shift-invariant subgroups of $\F_p^{\Z^d}$
(characteristic $p$, where $p\geq 2$ is a prime integer).

For these actions, mixing of all orders is not satisfied. Nevertheless, it is possible to show that non-mixing configurations are sparse
(Section \ref{special}). This was shown for a particular case of our model (Ledrappier's system) by D. Arenas-Carmona, D. Berend and V. Bergelson
in \cite{ArBeBe08}. We borrow from their paper, a source of inspiration for us, the term "almost mixing of all orders" used in the title.

The scarcity of non-mixing configurations allows to apply the cumulant method as in \cite{CohCo15} to prove the Central Limit Theorem
for different methods of summation (Section 3).

The last section (appendix) is devoted to reminders on algebraic endomorphisms of compact abelian groups.

\section{\bf Shift-invariant subgroups of $\F_p^{\Z^d}$ and a class of endomorphism} \label{shiftInv}

In this section, we recall some facts about shift-invariant subgroups of $\F_p^{\Z^d}$ (cf. \cite{Sch95}) and we define a class of endomorphisms
of these groups.

\subsection{\bf Shift-invariant subgroups of $\F_p^{\Z^d}$}

\

{\it Notations}

Let $p > 1$ be a prime integer fixed once for all and let $\F_p$ denote the finite field $\Z/p\Z$.
For all integers $a, b$, we have $a^p = a \mod p$, $(a+b)^p = a^p + b^p \mod p$.
Underlined symbols will represent vectors or tuples. The element $(0,0,...,0)$ is represented by $\0$.
For $d \geq 1$, if $J$ is a set of indices, a $|J|$-tuple of elements of $\Z^d$ or $\N^d$ ($\N$ includes 0) is written $\a_J = (\a_j, j \in J)$.
The coordinates are denoted by $a_{j,k}, \, j \in J, \, k = 1, ...,d$. The notation $\x_J$ represents the variable
$(x_j, j \in J)$ or the formal product $\prod_{j \in J} x_j$.

We denote by $G_0^{(d)}$, or simply $G_0$, the compact abelian group $\F_p^{\Z^d}$ (with coordinatewise addition and endowed with the product topology)
identified with the ring ${\Cal S}_d = \F_p[[x_1^\pm, ..., x_d^\pm]]$ of formal power series in $d$ variables with coefficients in $\F_p$.

An element $\zeta = (\zeta_\k, \,\k \in \Z^d)$ in $G_0$ is represented by the formal power series\footnote{\ We write $\x$ for $(x_1, \cdots, x_d)$
as well as for $x_1 \cdots x_d$, $\x^{-\k}$ for $x_1^{-k_1} \cdots \, x_d^{-k_d}$, $\zeta_\k$ or $c(\zeta, \k)$ for the coefficients of the series
$\zeta(\x)$.} with coefficients in $\F_p$: $\zeta(\x) = \sum_{\k \in \Z^d} \zeta_\k \, \x^{-\k}$.
For $j=1, ..., d$, the shift $\sigma_j$ on $\F_p^{\Z^d}$ corresponds to the multiplication by $x_j$: $\zeta(\x) \to x_j \zeta(\x)$.

{\it Polynomials with coefficients in $\F_p$ and characters}

The ring $\F_p [x_1^\pm, ..., x_d^\pm]$ of Laurent polynomials in $d$ variables with coefficients in $\F_p$ is denoted by ${\Cal P}_d$ .
For $d=1$, we write simply $\Cal P$. A Laurent polynomial $P \in {\Cal P}_d$ reads
$$P(x_1, ...,x_d) = \sum_{\k \in S(P)} c(P, \k) \, \x^\k,$$
where $S(P)$, called the support of $P$, is the finite set $\{\k: c(P, \k) \not = 0\}$.

For $P \in {\Cal P}_d$ and $\zeta \in {\Cal S}_d$, the product $P\zeta$ is well defined:
$$(P\zeta)(\x) = (\sum_{\el \in \Z^d} c(P, \el) \, \x^\el) \, (\sum_{\el' \in \Z^d} \zeta_{\el'} \, \x^{-\el'})  =
\sum_{\k \in \Z^d} \ \bigl (\sum_{\j \in S(P)} c(P, \j) \, \zeta_{\k + \j}\bigr) \, \x^{-\k}.$$

The dual $\widehat {G_0^{(d)}}$ of $G_0^{(d)}$ can be identified with ${\Cal P}_d$:
for any character $\chi$ on $\F_p^{\Z^d}$ there is a polynomial $P \in \Cal P_d$ such that
\begin{eqnarray}
\chi(\zeta) = \chi_P(\zeta) := e^{{2 \pi \over p} i \sum_{\el \in S(P)} c(P, \el) \, \zeta_\el}
= e^{{2 \pi \over p} i c(P\zeta, \0)}. \label{characp}
\end{eqnarray}

{\it Shift-invariant subgroups of $\F_p^{\Z^d}$}

Let $G \subset G_0$ be a shift-invariant closed subgroup of $\F_p^{\Z^d}$.
The annulator $G^\perp$ of $G$ in $\widehat G_0$ is $\{P: \chi_P(\zeta) = 1, \ \forall \zeta \in G\} =
\{P: c(P\zeta, \0) = 0, \ \forall \zeta \in G\}$. Since $G$ is shift-invariant, if $c(P\zeta, \0) = 0, \ \forall \zeta \in G$,
the same relation is satisfied for $\x^\k \, \zeta(\x)$, $\forall \k \in \Z^d$, which implies $P\zeta = 0, \ \forall \zeta \in G$.

Therefore $G^\perp$ can be identified with the ideal ${\Cal J} = \{P \in \Cal P_d: P\zeta = 0, \ \forall \zeta \in G\}$.
Since, by duality in $\F_p^{\Z^d}$, we have $G = (G^\perp)^\perp$ (see Appendix), this shows that $G = G_{\Cal J}$ where
\begin{eqnarray}
G_{\Cal J} = \{\zeta: P\zeta = 0, \ \forall P \in  \Cal J\}. \label{Gideal}
\end{eqnarray}
Conversely, for every ideal ${\Cal J} \subset \Cal P_d$, (\ref{Gideal}) defines a shift-invariant subgroup $G_\Cal J$ of $G_0$.

The dual of $G_{\Cal J}$ is isomorphic to the quotient $\widehat G_0 / G_{\Cal J}^\perp$, i.e., $\widehat G_{\Cal J}  = {\Cal P_d} / {\Cal J}$.

\vskip 3mm
\subsection{\bf Endomorphisms of $\F_p^{\Z}$ and their invertible extension} \label{actFp}

\

For $R \in \Cal P$, let $\gamma_R$ be the endomorphism of $K:=\F_p^{\Z}$ defined by $R$.
The action of $\gamma_R$ on characters is the multiplication $P \to RP$.
If $R \not = 0$, the surjectivity of $\gamma_R$ on $K$, or equivalently the injectivity of the action of $\gamma_R$ on the dual group $\hat K$,
is clear, since $RP \equiv 0$ if and only of $P \equiv 0$.

The invertible extension of $\gamma_R$
can be constructed by duality from the action $F \to RF$ on the ring $\cal F_R[x^\pm]$ of fractions of the form $F(x) = {P(x) \over R(x)^\ell}$,
$P \in \Cal P$, $\ell \in \N$. The invertible extension of $\gamma_R$ is the dual action on the compact group dual of the discrete additive group
$\cal F_R$.

An isomorphic version of the invertible extension is obtained in the following way.
Let us consider the subgroup $G_{\Cal J}$ of $G^{(2)}$, where $\cal J$ is the ideal in $\cal P_2$ generated by the polynomial $x_2 - R(x_1)$.
Observe that the homomorphism $h_R$ from $\Cal P_2$ to $\Cal F_R$ defined by $P(x_1, x_2) \to P(x_1, R(x_1))$ is surjective
and has for kernel the ideal $\cal J$ (Lemma \ref{basis0} below). Therefore we get an isomorphism between $\cal P_2 \mod \cal J$ and $\cal F_R$.
The shift $\sigma_2$ on the second coordinate is the invertible extension of $\gamma_R$.

{\it Multidimensional action}

This construction can be extended to a multidimensional action. We start with the group $K= \F_p^{\Z}$ and with $d$ polynomials $R_1, ..., R_d$
in $\F_p[x^{\pm}]$. We can add the polynomial $R_0(x) = x$ to the list.

The family $\gamma_{R_1}, ..., \gamma_{R_d}$ generate an $\N^d$-action on $K$ by algebraic endomorphisms.
The dual action on $\hat K= \cal P$ is $\el \to (P \to \uR^\el P)$, where $\uR^\el(x) := R_1^{\ell_1}(x)... R_d^{\ell_d}(x)$.

The natural extension is constructed as follows.
Let $\Cal F_{\Cal R}[x^\pm]$ be the ring of rational fractions in one variable with denominators $R_j(x)$, $j= 1, ..., d$, i.e.,
the discrete group of rational fractions with coefficients in $\F_p$ of the form
${V \over R_1^{\ell_1} \dots R_d^{\ell_d}}, \, V \in \Cal P[x^\pm], \, \ell_1, ..., \ell_d \in \N$.

Using duality, the natural extension $\tilde K$ of $K$ (with respect to the endomorphisms $\gamma_{R_1}, ..., \gamma_{R_d}$)
can be built as the dual of $\Cal F_{\Cal R}[x^\pm]$ view as an additive group.

As above, we can get an isomorphic version of the invertible extension in a shift-invariant subgroup of $G_0^{(d+1)}$,
namely the shift-invariant subgroups $G_{\Cal J}$ of $G_0^{(d+1)}$, where $\Cal J$ is the ideal generated in $\cal P_{d+1}$
by $x_{j+1} -R_j(x_1)$, $j = 1, ..., d$.

\begin{lem} \label{invertExt0} The shifts $\sigma_2, ..., \sigma_{d+1}$ are the invertible extensions of the endomorphisms
$\gamma_{R_1}, ..., \gamma_{R_d}$ acting on $\F_p^{\Z}$ and generate a $\Z^d$-action on $G^{(d+1)}_{\cal J}$.
\end{lem}
\proof \ If $\Gamma$ is in $\cal P_{d+1} = \F_p[x_1^\pm, x_2^\pm,, ..., x_{d+1}^\pm]$, let $h_\Cal R (\Gamma)$ be the rational fraction
\begin{eqnarray}
h_\Cal R (\Gamma)(x) = \Gamma(x, R_1(x),..., R_d(x)). \label{homomh}
\end{eqnarray}
The map $h_\Cal R$ is a surjective homomorphism from $\cal P_{d+1}$ to $\Cal F_{\Cal R}[x^\pm]$.
The homomorphism $\widetilde h$ defined by $\widetilde h \, (\Gamma \mod {\Cal J}) := h(\Gamma)$ is well defined,
since $\Gamma \in {\Cal J}$ implies $h(\Gamma)=0$. By Lemma \ref{basis0} below, it is an isomorphism between
$\Cal P_{d+1} \mod {\Cal J}$ and ${\Cal P}$.

The multiplication by $x_{j+1}$ on ${\Cal P}_{d+1} \mod {\Cal J}$ corresponds by $\tilde h$ to the multiplication by $R_j$
on ${\Cal P}$ and the $\Z^{d}$-action generated by the shifts $\sigma_{2}, ..., \sigma_{d+1}$ on
$G_{{\Cal J}}$ has the $\N^{d}$-action generated by $\gamma_{R_1}, ..., \gamma_{R_d}$ on $K$ as a factor through the map $\widetilde h$. \eop

\begin{lem} \label{basis0} The polynomials $L_j(\x) := x_{j+1} - R_j(x_1)$, $j= 1, ..., d$, form a basis of $\Ker h_\Cal R$.
\end{lem}
\proof \ Let us take for simplicity $d=3$.
If $P$ is in $\Ker h_\Cal R$, then $x_2 = R_2(x)$ is a root of the polynomial $Q_{x}(x_2)$ in $x_2$ defined by
$Q_{x}(x_2) = P(x, x_2, R_3(x)) - P(x, R_2(x), R_3(x))$. Therefore, there is $V_x$ such that
$P(x, x_2, R_3(x)) = P(x, R_2(x), R_3(x)) + V_x(x_2) (x_2 - R_2(x)) = V_x(x_2) (x_2 - R_2(x))$. The last equality follows from $P \in \Ker h_\Cal R$.

Now, $R_3(x)$ is a root of the polynomial $Q_{x, x_2}(x_3)$ in $x_3$ defined by $P(x, x_2, x_3) - P(x, x_2, R_3(x))$.
There is $W_{x, x_2}$ such that $P(x, x_2, x_3) = P(x, x_2, R_3(x)) + W_{x, x_2}(x_3) (x_3 - R_3(x))$. Put together, it gives:
$P(x, x_2, x_3) = V_x(x_2) (x_2 - R_2(x)) + W_{x, x_2}(x_3) (x_3 - R_3(x)).$

$V_x$ and $W_{x, x_2}$ can be written as polynomials, respectively $V(x,x_2)$ and $W(x, x_2, x_3)$. We have $P = V \, L_2 + W \, L_3$. \eop

Ledrappier's example (\cite{Le78}) corresponds to $p=2$, $d = 1$, $R_1(x) = 1+x$.
In this case, the invertible extension of $\gamma_R$ is given by the shift action on the second coordinate for the shift-invariant group
$G_{\Cal J}$ associated to the ideal ${\Cal J}$ generated by the polynomial $1+x_1+x_2$. The group $G_{\Cal J}$ is the set of configurations
$\zeta$ in $\F_2^{\Z^2}$ such that $\zeta_{n,m} + \zeta_{n+1,m} + \zeta_{n,m+1} = 0 \mod 2$, $\forall (n,m) \in \Z^2$.
The $\Z^2$-shift-action on $G_{\Cal J}$ endowed with its Haar measure is not $r$-mixing for $r \geq 3$, a fact which is general for the model
described above.

{\it Generalization: endomorphisms of $G_{\Cal J}$ and their invertible extension}

For an integer $d_1 \geq 1$, let $G_{\Cal J}$ be a shift-invariant subgroup of $G_0^{(d_1)}$.

Every polynomial $R$ in ${\Cal P}_d$
defines an endomorphism of $G_{\Cal J}$, $\gamma_R: \zeta(\x) \to R(\x) \, \zeta(\x)$.
Indeed, if $\zeta $ is such that $P\zeta = 0$, then $PR\zeta = R P\zeta = 0$.
 The dual action of $\gamma_R$ on $\hat G_{\Cal J}$ is the map $P \mod {\Cal J} \to RP \mod {\Cal J}$.

Any family $R_1, ..., R_{d_2}$, $d_2 \geq 1$,
of Laurent polynomials in $\x= (x_1, ..., x_{d_1})$ defines an $\N^{d_2}$-action by commuting endomorphisms $\gamma_{R_j}$ of $G_{\Cal J}$ .

The natural invertible extension of this action to a $\Z^{d_2}$-action by algebraic automorphisms of an extension of $G_{\Cal J}$ can be obtained
as above in the following way.

Let $x_{d_1+1}, ..., x_{d_1+d_2}$ be additional coordinates and consider $G_0^{(d_1+d_2)} = \F_p^{\Z^{d_1+d_2}}$.
The ideal ${\Cal J}'$ in ${\Cal P}_{d_1+d_2}$ generated by $\Cal J$ (embedded in $\Cal P_{d_1+d_2}$) and by the polynomials
$x_{d_1+1} - R_1(\x), ..., x_{d_1+d_2} - R_{d_2}(\x)$ defines a shift-invariant subgroup $G_{{\Cal J}'}$ of $G_0^{(d_1+d_2)}$.

Let us consider the surjective homomorphism $h$ from the ring $\Cal P_{d_1+d_2}$ of polynomials in $d_1+d_2$ variables
to the ring ${\Cal P}_{d_1}$ of polynomials in $d_1$ variables defined by $h(Q)(\x) = Q(\x, R_1(\x), ..., R_{d_2}(\x))$.

The homomorphism $\widetilde h$ defined by $\widetilde h \, (Q \mod {\Cal J}') := h(Q) \mod {\Cal J}$ is well defined,
since $Q \in {\Cal J}'$ implies $h(Q) \in {\Cal J}$. Using an extension of Lemma \ref{basis0} below, it can be shown
that it is an isomorphism between $\Cal P_{d_1+d_2} \mod {\Cal J}'$ and ${\Cal P}_{d_1} \mod {\Cal J}$.

The multiplication by $x_{d+j}$ on ${\Cal P}_{d_1+d_2} \mod {\Cal J}'$ corresponds by $\tilde h$ to the multiplication by $R_j$
on ${\Cal P}_{d_1} \mod {\Cal J}$. In other words, the $\Z^{d_2}$-action generated by the shifts $\sigma_{d_1+1}, ..., \sigma_{d_1+d_2}$ on
$G_{{\Cal J}'}$ has the $\N^{d_2}$-action generated by the endomorphisms $R_{d_1+1}, ..., R_{d_1+d_2}$ on $G_{{\Cal J}}$ as a factor through
$\widetilde h$.

The action of the shifts on $G_{{\Cal J}'}$ generate a $\Z^{d_1 + d_2}$-action, invertible extension of the action
generated on $G_{{\Cal J}}$ by multiplication by $x_1, ..., x_{d_1}, R_1(\x), ..., R_{d_2}(\x)$.

In the sequel we restrict the previous model to the case $\Cal J = \{0\}$. Moreover, although we think that the methods used below can be extended
to $d_1 >1$, we take  $d_1= 1$.
\vskip 2mm
{\it Total ergodicity}

Suppose that the polynomials $R_j, j=1, ..., d,$ are pairwise relatively prime of degree $\geq 1$. Then the family $(R_j, j=1, ..., d)$ generates
an $\Z^d$-action on $K= \F_p^{\Z}$ by endomorphisms, which extends to a $\Z^d$-action $(A^\el, \el \in \Z^d)$ on the natural
extension $\tilde K$ of $K$ which is {\it totally ergodic} (i.e. such that $A^\el$ on $(\tilde K, \tilde \mu)$ is ergodic for every $\el \in \Z^d \stm0$).

{\it Example}: (with $d = 3$) We take $p=2$, $R_0(x) = x$, $R_1(x) = 1+x$, $R_2(x) = 1+x+x^2$.
The orbits on the set of non trivial characters of the generated $\Z^3$-action are infinite by primality of the polynomials $x$, $1+x$, $1+x+x^2$.
Therefore, we get a 2-mixing $\Z^3$-action, hence a $\Z^3$-action with Lebesgue spectrum on $L_0^2(\tilde \mu)$,
where $\tilde \mu$ is the Haar measure on $\tilde K$.

\subsection {\bf Non $r$-mixing tuples}

\

Let us briefly recall the relation between $r$-mixing for an action by algebraic endomorphims and $S$-unit equations.
A general measure preserving $\N^d$-action $(T^\el)_{\el \in \N^d}$ on a probability measure space $(X, \mu)$
is mixing of order $r \geq 2$ if, for any $r$-tuple of bounded measurable functions $f_1, ..., f_r$ on $X$ with 0 integral
and for every $\varepsilon >0$, there is $M \geq 1$ such that
\begin{eqnarray}
\|\el_j - \el_{j'}\| \geq M, \forall j \not = j' \, \Rightarrow |\int T^{\el_1} f_1 ... T^{\el_r} f_r \, d\mu| < \varepsilon. \label{nonMix}
\end{eqnarray}
When $(X, \mu)$ is a compact abelian group $G$ with its Haar measure, one easily checks by approximation that
mixing of order $r$ for an $\N^d$-action generated by algebraic endomorphisms $T_1, ..., T_d$ is equivalent to:
{\it for every set ${\Cal K} = \{\chi_1,..., \chi_r\}$ of $r$ characters different from the trivial character $\chi_0$,
there is $M \geq 1$ such that $\|\el_j - \el_{j'}\| \geq M$ for $j \not = j'$ implies $T^{\el_1} \chi_1... T^{\el_r} \chi_r \not = \chi_0$.}

The "non-mixing" $r$-tuples in $\N^d$ for $\Cal K$ are the $r$-tuples in the set
\begin{eqnarray}
\Phi({\Cal K}, r) := \{(\el_1, ..., \el_r): \, T^{\el_1} \chi_1  ... T^{\el_r} \chi_r = \chi_0\}. \label{nonMix1}
\end{eqnarray}

\vskip 2mm
\goodbreak
{\it Example: action by $\times 2$, $\times 3$ on $\T^1$}

Let us illustrate the question of mixing on an example in the connected case: the action $\times 2$, $\times 3$ on $\T^1$.
A set ${\Cal K}$ of non zero characters on the torus is given by an $r$-tuple $\{k_1,..., k_r\}$ of non zero integers.
By putting $\el_j = (a_j, b_j)$, Equation (\ref{nonMix1}) for the action by 2 and 3 reads $k_1 2^{a_1}3^{b_1} + ...+ k_r2^{a_r}3^{b_r}=0$,
which leads to consider equations of the form:
\begin{eqnarray} k_1 2^{a_1}3^{b_1} + ...+ k_r2^{a_r}3^{b_r}=1, \ ((a_1,  b_1), ..., (a_r, b_r)) \in (\Z^2)^r. \label{2and3} \end{eqnarray}
It is known that, for a given set $\{k_1, ..., k_r\}$, the number of $r$-tuples $((a_1,  b_1), ..., (a_r, b_r))$ solutions of (\ref{2and3}),
such that no proper subsum vanishes, is finite (cf. Theorem \ref{thm1.1}). It implies that the
$\Z^2$-action generated by the invertible extension $\times 2, \times 3$ is mixing of all orders.
This mixing result is a special case of a general theorem of K. Schmidt and T. Ward (1992):
\begin{thm} (\cite{SchWar93}) \label{SchmWrd} Every 2-mixing $\Z^d$-action by automorphisms on a compact connected abelian group G is
mixing of all orders.
\end{thm}
The proof of Theorem \ref{SchmWrd} relies on a result on $S$-unit equations (Schlickewei (1990)).
Let us mention the following version of results on $S$-unit equations in characteristic $0$:

Let $\F$ be an algebraically closed field of characteristic 0, $\F^*$ its multiplicative group of nonzero elements.
Let $\Gamma$ be a subgroup of $(\F^*)^r$.
\begin{thm} \label{thm1.1} (J.-H. Evertse, J.-H. Schlickewei, W. M. Schmidt \cite{EvScSch02}) If the rank of $\Gamma$ is finite,
for $(k_1,..., k_r) \in (\F^*)^r$, the number of solutions $(\gamma_1, ..., \gamma_r) \in \Gamma$ of equation
\begin{eqnarray}
k_1 \gamma_1 + ... +k_r \gamma_r = 1, \label{(1.1)}
\end{eqnarray}
such that $\sum_{i \in I}  k_i \gamma_i \not = 0$ for every nonempty subset $I$ of $\{1,..., r\}$, is finite.
\end{thm}

{\it "Non-mixing" $r$-tuples for the action of $R_1, ..., R_d$}

The situation is different in characteristic $\not = 0$, where there can exist infinitely many solutions for equations of the type (\ref{(1.1)}).
In the non connected case (for example for endomorphisms of shift-invariant subgroups of $\F_p^{\Z^d}$),
this implies the existence of infinitely many non-mixing $r$-tuples, for $r \geq 3$.

Our goal is to show that, however, these non-mixing $r$-tuples
for the $\N^d$-actions described in Subsection \ref{actFp} above are rare in a sense (hence these actions are "almost mixing of all order")
(cf. \cite{ArBeBe08} and D. Masser's works about non-mixing $r$-tuples).

Our framework is the setting introduced previously. We consider the $\N^d$-action on $\F_p^\Z$ defined by ${\Cal R} = (R_1, ..., R_d)$.
A finite set of characters is given by a finite family of polynomials $P_1, ..., P_r$. For such a set, a non-mixing $r$-tuple
of the action is an $r$-tuple $(\a_1, ..., \a_r) \in (\N^d)^r$ such that in $\F_p[x^\pm]$
\begin{eqnarray}
P_1(x) \prod_{i=1}^d R_i(x)^{a_{1,i}} + ... + P_r(x) \prod_{i=1}^d R_i(x)^{a_{r,i}} = 0. \label{equaBad}
\end{eqnarray}
Equation (\ref{equaBad}) is analogous to the previous $S$-unit equation (\ref{(1.1)}), but in characteridtic $p \not = 0$.
Observe that, for a given family ($P_j)$, the equation can be reduced to the case where the $P_j$'s are scalars: it suffices to enlarge the family $\Cal R$
by adding the $P_j$'s to $\Cal R$.

Replacing $R_j$ in $\F_p[x^\pm]$ by $\widetilde R_j$, the polynomial in $\F_p[x]$ such that $\widetilde R_j(x) = x^\ell R_j(x)$
where $\ell \geq 0$ is minimal, we can also suppose that the polynomials $R_j$ are in $\F_p[x]$.

To count non-mixing $r$-tuples (for the action of $\Cal R$ on $K$ or of the shifts on the natural invertible extension), in the next section we will
study polynomials $\Gamma$ which belong to $\Ker(h_{\Cal R})$ where $h_{\Cal R}$ is the homomorphism defined by (\ref{homomh}).

\section{\bf Basic special $\Cal D$-polynomials} \label{special}

\subsection{\bf Decomposition of special $\Cal D$-polynomials}

\subsubsection{\bf Preliminary notations and results}

\

In this section, we extend results shown for Ledrappier's example of \cite{ArBeBe08} to the general model introduced in the first section.
We start the proof of the main theorem (Theorem \ref{sumBasic0}) with some notations and preliminary results.

{\it Notations}: We denote by $\Upsilon$ the set of all monic (i.e., with leading coefficient equal to 1) prime polynomials in one variable over $\F_p$.

For $U \in \Cal P$, $\Upsilon(U)$ denotes the set of its prime monic factors.
If $U$ is a constant $\not = 0$, we set $\Upsilon(U) = \{1\}$.
If $\Cal S$ is a family of polynomials in one variable,
$\Upsilon(\Cal S) := \bigcup_{U \in \Cal S} \Upsilon(U)$ is the set of their prime monic factors.

We denote by ${\Cal Q}_0$ the ring of (Laurent) polynomials, with coefficients in $\F_p$, in the  variables $x_\rho$ indexed by $\rho \in \Upsilon$.
By definition, for every $\Gamma$ in ${\Cal Q}_0$, there is a finite subset $J(\Gamma)$ of $\Upsilon $ such that
$\Gamma$ is a polynomial in the variables $x_\rho, \rho \in J(\Gamma)$, and reads (in reduced form):
\begin{eqnarray}
\Gamma(\x) &&= \sum_{\a \in \Z^{J(\Gamma)}} d(\a) \, \prod_{\rho \in J(\Gamma)} \, x_\rho^{a_\rho}, \text{ with } d(\a) \in \F_p. \label{Qexpr0}
\end{eqnarray}
The term ``reduced" means that a product $\prod_{\rho \in J(\Gamma)} \, x_\rho^{a_\rho}$ in the above formula appears
only once for a given $\a \in \N^{J(\Gamma)}$ with a coefficient $d(\a) \not = 0$, except for the 0 polynomial.
Most of the time it will be enough to consider polynomials with non negative exponents.

If a polynomial $\Gamma$ is expressed in a non reduced form, its expression in reduced form (possibly the 0 polynomial) is
$$red(\Gamma)(\x) = \sum_{\b \in \Z^{J(\Gamma)}} (\sum_{\a: \, \a=\b}\ d(\a)) \, \prod_{\rho \in J(\Gamma)} \, x_\rho^{b_\rho}.$$

An element $(\alpha_1, ..., \alpha_d)$ of $\{0, ..., p-1\}^d$ (identified to $\F_p^d$) is denoted by $\al$.
For $\Gamma$ given by (\ref{Qexpr0}), we call $\al$-homogeneous component of $\Gamma$, for $\al \in \F_p^{J(\Gamma)}$
the sum:
\begin{eqnarray}
\Gamma_{\al}(\x) = \sum_{\b \in \Z^{J(\Gamma)}} d(p \, \b + \al) \, \prod_{\rho \in J(\Gamma)} \, x_\rho^{p \, b_\rho+\alpha_\rho}.  \label{Qexpr}
\end{eqnarray}
There is a homomorphism $h: \Gamma \to h(\Gamma)$, denoted also $\Gamma \to \widehat \Gamma$, from ${\Cal Q}_0$ to $\Cal P$, defined by
\begin{eqnarray}
\Gamma(\x) = \sum_{\a \in \Z^{J(\Gamma)}} d(\a) \, \prod_{\rho \in J(\Gamma)} \, x_\rho^{a_\rho}
\ \to \ h(\Gamma)(x) := \sum_{\a \in \Z^{J(\Gamma)}} d(\a) \, \prod_{\rho \in J(\Gamma)} \, \rho(x)^{a_\rho}. \label{homh1}
\end{eqnarray}

We consider also the ring ${\Cal Q}_1$ of polynomials $\Gamma$ in the variables $x_\rho, \rho \in \Upsilon$, with coefficients
in $\F_p[x]$:
\begin{eqnarray}
&&\Gamma(x, \x) = \sum_{\a \in \Z^{J(\Gamma)}} d(\a) \, U_\a(\x) \, \prod_{\rho \in J(\Gamma)}
\, x_\rho^{a_\rho}. \label{Q1}
\end{eqnarray}

{\bf Definitions}: If $\Cal D$ is a finite subset of $\Upsilon$ (i.e., a finite set of prime polynomials), a polynomial in ${\Cal Q}_1$ of the form
\begin{eqnarray}
\Gamma(x, \x) &&= \sum_{\a \in \Z^{\Cal D}} d(\a) \, U_\a(x) \, \prod_{\rho \in \Cal D} \, x_\rho^{a_\rho}, \text{ with } d(\a) \in \F_p,
\ U_\a(x) \text{ monic}, \label{QU1}
\end{eqnarray}
is called a {\it $\Cal D$-polynomial}. It is called a {\it special $\Cal D$-polynomial} if it satisfies
\begin{eqnarray}
\sum_{\a \in \Z^{\Cal D}} d(\a) \, U_\a(x) \, \prod_{\rho \in \Cal D} \, \rho(x)^{a_\rho} = 0. \label{spQU1}
\end{eqnarray}
For $U \in \Cal P$, if $U(x) = c(U) \prod_{\rho \in \Upsilon(U)} \, \rho(x)^{\theta_\rho(U)}, \ c(U) \in \F_p$,
is the factorization of $U$ into prime monic factors, we put
$$\Psi(U)(\x) = c(U) \prod_{\rho \in \Upsilon(U)} x_\rho^{\theta_\rho(U)}.$$
For example, for $p=2$ and $U(x) = x^3 + x^5$, denoting by $\rho_1, \rho_2$ the polynomials $x$ and $1+x$,
we get $\Psi(U)(\x) = x_{\rho_1}^3x_{\rho_2}^2$.

Observe that $\Psi(U)(\x) - U(x)$ is a special $\Upsilon(U)$-polynomial.

We define now a map $\Gamma \to \Psi(\Gamma)$, also denoted $\Gamma \to \widetilde \Gamma$, from ${\Cal Q}_1$ to ${\Cal Q}_0$,
which maps $\Gamma$ given by (\ref{Q1}) to the (not necessarily reduced) polynomial $\widetilde \Gamma$:
\begin{eqnarray}
&&\Psi(\Gamma)(\x) = \widetilde \Gamma(\x) = \sum_{\a \in \Z^{J(\Gamma)}} d(\a) \,\Psi(U_\a)(\x) \, \prod_{\rho \in J(\Gamma)}
\, x_\rho^{a_\rho}. \label{psimap}
\end{eqnarray}
If $\Gamma$ is a special $\Cal D$-polynomial, then $\widetilde \Gamma$ is a special $\Cal D \cup_\a \Upsilon(U_\a)$-polynomial.
Denoting by $r(\Gamma)$ the number of terms of $\Gamma$ and $S(\Gamma)$ its support,
observe that $\widetilde \Gamma - \Gamma$ is a sum of $r(\Gamma)$ special $\bigcup_\a \Upsilon(U_\a)$-polynomials:
\begin{eqnarray}
\widetilde \Gamma(\x) - \Gamma(\x) = \sum_{\a \in S(\Gamma)} c(\a) \, [\prod_{\rho \in \Upsilon(U_\a)}x_\rho^{\theta_\rho(U_\a)} - U_\a(x)]\,
\prod_{\rho \in J(\Gamma)} \, x_\rho^{a_\rho}. \label{basic0}
\end{eqnarray}

{\it Basic special $\Cal D$-polynomials}

Let $\Cal D$ be any family of prime polynomials containing the polynomial $x \to x$. The polynomials $x_\rho - \rho(x), \rho \in \Cal D$,
are called {\it basic special $\Cal D$-polynomials} (abbreviated in "bs~$\Cal D$-polynomial").
We say that a polynomial $\Gamma$ is {\it shifted} from $\Gamma_0$ if $\Gamma(\x) = \x^\a \, \Gamma_0(\x)$
for some monomial $\x^\a$. We will use the following elementary lemma:
\begin{lem} \label{sumverybasic} For any monic polynomial $U$ in one variable, $\Psi(U)(\x) - U(x)$ is a sum of polynomials shifted from basic special
$\Upsilon(U)$-polynomials.
\end{lem}
\proof If $U$ is a power of a prime polynomial, $U(x) = \rho(x)^b$, $b \geq 1$, then we use:
$$x_\rho^b - \rho(x)^b = (\sum_{k=0}^{b-1} x_\rho^{b-k-1} \, \rho(x)^k) \, (x_\rho - \rho(x)).$$
The general case follows from the formula $Y^b Z^c - y^b z^c = (Y^b - y^b) \, Z^c + y^b(Z^c - z^c)$ by induction. \eop

A polynomial $\Lambda$ is called {\it generalized basic special} $\Cal D$-polynomial (abbreviated in ``gbs~$\Cal D$-polynomial"),
if it is obtained from a basic special $\Cal D$-polynomial $\Delta$
by shift and dilation (exponentiation with a power of $p$ as exponent).

Therefore $\Lambda$ is a gbs~$\Cal D$-polynomial if there are
$\a \in \Z^d$, $t \geq 0$ and a bs~$\Cal D$-polynomial $\Delta = x_\rho - \rho(x)$ such that:
$$\Lambda(\x) = \x^\a \, (\Delta(\x))^{p^t} =  \x^{\a} \, (x_\rho^{p^t} + (- \rho(x))^{p^t}.$$

In the sequel, $\Cal R =(R_j, j= 1, ..., d)$ will be a {\it fixed finite family} of $d \geq 2$ distinct prime polynomials in one variable over $\F_p$.
If the polynomial $x \to x$ is not included in the family $\cal R$, we add it to the list.

For this fixed family, it is convenient to introduce another notation for polynomials in ${\Cal Q}_0$ depending
on the variables $x_\rho \in \Cal R$. We write them as polynomials in $d$ variables $x_i$:
\begin{eqnarray}
\Gamma(\x) &=& \sum_{\a \in \N^d} \, d(\a) \, \prod_{i=1}^d x_i^{a_i}. \label{gamma0}
\end{eqnarray}
The variable $x_i$ corresponds to the polynomial $R_i$. We will use the equivalent notations
$\x^\a$, $\prod_{i=1}^d x_i^{a_i}$ or $\prod_{\rho \in \Cal R} x_\rho^{a_\rho}$ (here $\Upsilon(\Cal R) = \Cal R$,
since the $R_j$'s are prime polynomials).

$\Gamma$ in (\ref{gamma0}) is written in its reduced form (a product $\prod_{i=1}^d x_i^{a_i}$
appears only once for a given $\a \in \N^d$). As above, $\Gamma$ reads as a sum of {\it $\al$-homogeneous components}:
\begin{eqnarray}
\Gamma(\x) &=& \sum_{\al \in \F_p^d} \, \Gamma_{\al}(\x) = \sum_{\al \in \F_p ^d} \ [\sum_{\b \in \N^d} c_{\b, \al} \,
\prod_{i=1}^d x_i^{p \, b_i+\alpha_i}] = \sum_{\al \in \F_p ^d} \, (\prod_{i=1}^d x_i^{\alpha_i}) \, \overline \Gamma_{\al}(\x), \label{homogen0}\\
&& \text{ with }
\overline \Gamma_\al(\x) := (\prod_{i=1}^d x_i^{-\alpha_i}) \, \Gamma_{\al}(\x) = \sum_{\b \in \N^d} \, c_{\b, \al} \prod_{i=1}^d x_i^{p \, b_i}.
\label{overP}
\end{eqnarray}
We denote by $r(\Gamma_\al)$ the number of monomials in the sum $\Gamma_{\al}$. The {\it length} of $\Gamma$ is the number $r(\Gamma)$ of its monomials.
It is the cardinal of the support of $\Gamma$.

\vskip 3mm
{\it The map $\Gamma \to \check \Gamma$}

In case the $R_i$'s are monic polynomials non necessarily prime, we use the reduction to the prime case given by the following map.
Let $R_i = \prod_{\rho \in \Upsilon(\Cal R)} \rho^{b_{i, \rho}}$. The map  $\Gamma \to \check \Gamma$ is defined by
\begin{eqnarray}
\Gamma(\x) = \sum_{\a \in S(\Gamma)} c(\a) \, \prod_{i=1}^d x_i^{a_i} \to \check \Gamma(\x_{\Upsilon(\Cal R)}) = red \, (\sum_{\a \in S(\Gamma)} c(\a)
\, \prod_{\rho \in \Upsilon(\Cal R)} x_\rho^{\sum_{i=1}^d a_i b_{i, \rho}}).
\end{eqnarray}
If $\Gamma$ is such that $h_{\Cal R} (\Gamma) = 0$, i.e., $\sum_{\a \in S(\Gamma)} c(\a) \, \prod R_i(x)^{a_i} = 0$,
then $\check \Gamma$ is a special $\Upsilon(\Cal R)$-polynomial.

{\it The goal of this section is the study of the set of special $\Cal R$-polynomials.
Theorem \ref{sumBasic0} will show that, for every family $\Cal R$ of polynomials and every $r$, there is a finite constant
$t(r, \Cal R)$ and a finite family $\Cal E$ of polynomials in one variable containing $\Cal R$ such that every special $\Cal R$-polynomial
$\Gamma$ of length $r$ is a sum of at most $t(r,\Cal R)$ gbs~$\Cal E$-polynomials. The constant $t(r,\Cal R)$ does not depend on the degree
of the polynomial $\Gamma$.}

Let us now recall or mention some facts about polynomials over $\F_p$.
\begin{lem} \label{nonZeroDer} a) For any polynomials $A, B$, we have $(AB^p)' = A' \, B^p$.

b) A product of pairwise relatively prime polynomials is a $p$-th power if and only if each factor is a $p$-th power.

c) If $P$ is a (reduced) polynomial in one variable, then $P' = 0$ if and only if $P = U^p$ for some polynomial $U$.

d) If $V_1, ..., V_n$ are pairwise relatively prime polynomials which are not $p$-th powers, then $(\prod_{i=1}^n V_i)' \not = 0$.
\end{lem}
\proof {\it a)}, {\it b)} are clear. For {\it c)}, suppose that $P' = 0$, with $P(x) = \sum_k \sum_{\ell = 0}^{p-1} c(k, \ell) \, x^{pk+\ell}$,
then $0 = P'(x) = \sum_k \sum_{\ell = 1}^{p-1} \ell c(k, \ell) \, x^{pk + \ell}$ hence $P(x) = [\sum_k c(k, 0) \, x^{k}]^p$.

For {\it d)}, observe that  $(\prod_{i=1}^n V_i)' = 0$ implies that $\prod_{i=1}^n V_i$ is equal to $U^p$ for some polynomial $U$ by {\it c)},
which is impossible by the hypotheses on the $V_j$'s and {\it b)}. \eop

\vskip 3mm
\subsubsection{\bf Decomposition of special $\Cal R$-polynomials}

\

Let $\Gamma$ be a polynomial as in (\ref{homogen0}). With the notation (\ref{overP}), for $\be \in \F_p^d$ we put
\begin{eqnarray}
A_{\be}(\Gamma)(\x) &:=& \sum_\al \, (\prod_{i=1}^d R_i^{\alpha_i + \beta_i}(x))' \, \overline \Gamma_\al(\x), \label{Abeta} \\
B_\be(\Gamma)(\x) &:=& - (\prod_{i} R_i^{\beta_i}(x)) \, A_\0(\Gamma)(\x) =
- (\prod_{i} R_i^{\beta_i}(x)) \, \sum_\al \, (\prod_{i=1}^d R_i^{\alpha_i}(x))' \, \overline \Gamma_\al(\x), \label{B0P}\\
\Pi_\be(\Gamma)(\x) &:=&(\prod_{i} R_i^{\beta_i}(x))' \, \sum_{\al} \, (\prod_{i=1}^d R_i^{\alpha_i}(x)) \, \overline \Gamma_\al(\x). \label{PiP}
\end{eqnarray}
We assume that $\Gamma$ is a special $\Cal R$-polynomial, i.e., $\widehat \Gamma = 0$.

It follows that $A_{\be}(\Gamma)$ (hence also $B_\be(\Gamma))$
is a special $\Cal R$-polynomial. Indeed we have by Lemma \ref{nonZeroDer} a):
\begin{eqnarray*}
\widehat A_{\be}(\Gamma)&&= \sum_{\al} \, \sum_\a c(\b, \al) \,  (\prod_{i} R_i^{\alpha_i + \beta_i})' \, \prod_i R_i^{pb_i} \\
&&= [\sum_{\al} \, \sum_\a c(\b, \al) \,  (\prod_{i} R_i^{\alpha_i + \beta_i}) \, \prod_i R_i^{pb_i}]'
= [(\prod_{i} R_i^{\beta_i}) \, \widehat \Gamma]'= 0.
\end{eqnarray*}
From the identity $(\prod_i R_i^{\beta_i +\alpha_i})' =  (\prod_i R_i^{\beta_i})' \, (\prod_i R_i^{\alpha_i})
+ (\prod_i R_i^{\beta_i}) \, (\prod_i R_i^{\alpha_i})'$, we get $$\Pi_\be(\Gamma) = A_{\be}(\Gamma) + B_\be(\Gamma).$$

{\it Notation:} \ For a finite family of prime polynomials $\Cal D = \{S_i, i \in I(\Cal D)\}$ and $\be = (\beta_1, ..., \beta_{I(\Cal D)})$,
we put
\begin{eqnarray}
&&D_{\Cal D, \be, 0} := \prod_{i \in I(\Cal D)} S_i^{\beta_i},
\ D_{\Cal D, \be, 1} := (\prod_{i \in I(\Cal D)} S_i^{\beta_i})',  \label{Dbetak} \\
&&\zeta(\Cal D) := \Cal D \, \cup \, \bigcup_{\be \in \F_p^{I(\Cal D)}} \Upsilon(D_{\Cal D, \be, 1}). \label{SDmap}
\end{eqnarray}
If we iterate $k$-times the map $\zeta: \Cal D \to \zeta(\Cal D)$ starting from a finite family of
prime polynomials $\Cal R$, we get a finite family of prime polynomials denoted by $\zeta^k(\Cal R)$.

Remark that, if the derivatives of order 1 of products of polynomials in a family of prime polynomials $\Cal R$ do not contain prime factors
$\not \in \Cal R$, then $\zeta(\Cal R) = \Cal R$. This the case in few examples like for $p = 2$: $\Cal R = \{x, 1+x\}$ (Ledrappier's example),
$\Cal R = \{x, 1+x, 1+x+x^2\}$.

The map $\Psi$ (also denoted by $\widetilde \, $) defined in (\ref{psimap}) gives for $\Pi_\be(\Gamma), A_{\be}(\Gamma), B_\be(\Gamma)$:
\begin{eqnarray}
&&\Pi_\be(\Gamma)(\x) = D_{\Cal R, \be, 1}(x) \, \sum_{\al} \, (\prod_{i} R_i^{\alpha_i}(x)) \, \overline \Gamma_\al(\x) \  \to \nonumber \\
&&\widetilde \Pi_\be(\Gamma)(\x) = (\prod_{\rho \in \Upsilon(D_{\Cal R, \be,1})} x_\rho^{\theta_\rho(D_{\Cal R, \beta, 1})}) \,
\sum_{\al} \, (\prod_{i} x_i^{\alpha_i}) \, \overline \Gamma_\al(\x)
 = \Psi(D_{\Cal R, \be, 1})(\x) \, \sum_{\al} \, \Gamma_\al(\x), \label{TPiP} \\
&&A_{\be}(\Gamma)(\x) = \sum_{\al} \, D_{\Cal R, \be+\al, 1}(x)  \, \overline \Gamma_\al(\x) \ \to  \nonumber \\
&&\widetilde A_{\be}(\Gamma)(\x) = \sum_{\al} \, (\prod_{\rho \in \Upsilon(\Cal R, D_{\be+\al, 1})}
x_\rho^{\theta_\rho(D_{\Cal R, \be+\al, 1})}) \, \overline \Gamma_\al(\x)
= \sum_{\al} \, \Psi(D_{\Cal R, \be+\al, 1})(\x)  \, \overline \Gamma_\al(\x), \label{TAbeta} \\
&&B_\be(\Gamma)(\x)= D_{\Cal R, \be, 0}(x) \, A_\0(\x) =D_{\Cal R, \be, 0}(x) \, \sum_{\al} \, D_{\Cal R, \al, 1}(x)  \, \overline \Gamma_\al(\x)
\ \to \nonumber \\
&&\widetilde B_\be(\Gamma)(\x)
= \Psi(D_{\Cal R, \be, 0})(\x) \sum_{\al} \, \Psi(\Cal R, D_{\Cal R, \al, 1})(\x) \, \overline \Gamma_\al(\x). \label{TB0P}
\end{eqnarray}
The polynomials $\widetilde A_{\be}(\Gamma), \widetilde B_\be(\Gamma)$ are special $\zeta(\Cal R)$-polynomials (with more variables
than $\Gamma$ in general). This follows from (\ref{basic0}) and from the fact that $A_{\be}(\Gamma), B_\be(\Gamma)$ are special
$\Cal R$-polynomials, as was shown above.

{\it Reduction of the number of terms}

For $\gamma \in \F_p^d$, we define  $u(\gamma)$ by $u(\gamma)_i = 0 \text{ if } \gamma_i = 0, u(\gamma)_i = p - \gamma_i
\text{ if } \gamma_i = 1, ..., p-1$, $i=1, ..., d.$ We have:
$D_{\Cal R, \0, 1}(x) = 0, \ D_{\Cal R, \be_1 + \al, 1}(x) = 0, \text{ for } \al= u(\be_1)$.

If $\Gamma$ is not reduced to 0, by shifting $\Gamma$ by a monomial, we can assume that $\Gamma_{\0} \not = 0$. If $\Gamma$ does not reduce
to the single component $\Gamma_\0$, there is $\be_1 \not = \0 $ such that $\Gamma_{\be_1} \not = 0$.

If $\Gamma$ does not reduce to a single homogeneous component, we can optimize the choices of components in the decomposition
(see the proof of Theorem \ref{sumBasic0}). There are at most $p^d$ non zero homogeneous components.  We get
\begin{eqnarray}
&&r(\widetilde A_\0(\Gamma)) \leq (1 - \lambda_r)\,r, \ r(\widetilde A_{\be_1'}) \leq (1 - \mu_r) \, r,\\
&&\text{ with } \lambda_r = \max(p^{-d}, r^{-1}), \ \mu_r = \max((1- \lambda_r) p^{-d}, r^{-1}). \label{nbterm}
\end{eqnarray}

Suppose that $\Gamma_\0, \Gamma_{\be_1} \not = 0$. Let $\be= u(\be_1)$. The polynomials $\widetilde A_{\be}, \widetilde B_\be(\Gamma)$,
are special $\zeta(\Cal R)$-polynomials with strictly less terms than $\Gamma$.

For a family $\Cal R$, we get from the differences $\widetilde \Pi_\be(\Gamma) - \Pi_\be(\Gamma)$,
$\widetilde A_\be(\Gamma)  - A_\be$, $\widetilde B_{\be}(\Gamma) - B_{\be}(\Gamma)$ respectively the following special polynomial:
\begin{eqnarray}
&&\Delta_{\Cal R, \al, 0}^\be(\x) := \Psi(D_{\Cal R, \be, 1})(\x) \, \prod_{i=1}^d x_i^{\alpha_i}
- D_{\Cal R, \be, 1}(x)\, \prod_{i=1}^d R_i^{\alpha_i}(x), \, \al \in \F^d, \label{Delt1}\\
&&\Delta_{\Cal R, \be+\al, 1}(\x) := \Psi(D_{\Cal R, \be+\al, 1})(\x) -
D_{\Cal R, \be+\al, 1}(x), \, \al \in \F^d, \al \not = \be', \label{Delt2} \\
&&\Delta_{\Cal R, \al, k}^\be(\x) := \Psi(D_{\Cal R, \be, 0})(\x) \, \Psi(D_{\Cal R, \al, 1})(\x)
- D_{\Cal R, \be, 0}(x) \, D_{\Cal R, \al, 1}(x), \, \al \in \F^d \setminus \{\0\}. \label{Delt3}
\end{eqnarray}
Some polynomials in the list can be 0 and there can be redundancy. With the notation used in (\ref{nbterm}), the number of these polynomials is
$$\leq \sum_\al r(\Gamma_\al) + \sum_{\al \not = \be'} r(\Gamma_\al)  + \sum_{\al \not = \0} r(\Gamma_\al) = r + r(1-\mu_r) + r(1-\lambda_r)
\leq 3 \, r(\Gamma).$$
By Lemma \ref{sumverybasic} each of them can be expressed as a sum of shifted basic polynomials, with a number of terms bounded by a constant $C$.
They are then shifted by the corresponding $\overline \Gamma_\al$ associated to the $\al$-homogeneous component of $\Gamma$.

The results of these preliminaries are summarized in the following lemma:
\begin{lem} \label{oneStep} Let $\Gamma$ be a special $\Cal R$-polynomial of length $r$.

Let $A_{\be}(\Gamma), \Pi_\be(\Gamma), B_\be(\Gamma), \widetilde \Pi_\be(\Gamma), \widetilde A_{\be}(\Gamma), \widetilde B_\be(\Gamma)$
be defined respectively by (\ref{Abeta}), (\ref{PiP}), (\ref{B0P}), (\ref{TPiP}), (\ref{TAbeta}), (\ref{TB0P}). Then we have
\begin{eqnarray}
&&\Psi(D_{\Cal R, \be, 1}) \, \Gamma = \widetilde \Pi_\be(\Gamma) \label{decGamm1}\\
&&\quad \quad \quad = \widetilde A_{\be}(\Gamma) + \widetilde B_\be(\Gamma) + \widetilde \Pi_\be(\Gamma) - \Pi_\be(\Gamma)
+ A_{\be}(\Gamma) - \widetilde A_{\be}(\Gamma) + B_\be(\Gamma) - \widetilde B_\be(\Gamma). \label{decGamm2}
\end{eqnarray}

$\widetilde A_{\be}(\Gamma)$ and $\widetilde B_\be(\Gamma)$ are special
$\zeta(\Cal R)$-polynomials with a number of terms strictly less than the number of terms of $\Gamma$.

The differences $\widetilde \Pi_\be(\Gamma) - \Pi_\be(\Gamma)$, $A_{\be}(\Gamma) - \widetilde A_{\be}(\Gamma)$, $B_\be(\Gamma) - \widetilde B_\be(\Gamma)$
are sums of at most $C r$ gbs~$\zeta(\Cal R)$-polynomials.

The polynomial $\widetilde \Pi_\be(\Gamma)$ differs from $\widetilde A_{\be}(\Gamma) + \widetilde B_\be(\Gamma)$ by at most $3 C\, r(\Gamma)$
gbs~$\zeta(\Cal R)$-polynomials given by the decomposition of (\ref{Delt1}), (\ref{Delt2}), (\ref{Delt3}).
\end{lem}

Now we prove the main result of this section, which will be used to show that the non-mixing configurations are sparse for the actions that we consider.
\begin{thm} \label{sumBasic0} Let $r$ be an integer $\geq 2$. For every family $\Cal R =(R_j, j= 1, ..., d)$ of $d \geq 1$ polynomials,
there is a finite constant $t(r, \Cal R)$ and a finite family $\Cal E$ of polynomials in one variable containing $\Cal R$ such that every special
$\Cal R$-polynomial of length $ \leq r$ is a sum of at most $t(r,\Cal R)$ gbs~$\Cal E$-polynomials.

Moreover, $\Cal E = \zeta^{r_1}(\Cal R)$ for some $r_1 \leq r$ and there are two constants $K > 0$, $\theta \geq 2$, such that
$t(r, \Cal R) \leq K r^\theta$.
\end{thm}
\proof \ Let $H(r_0)$ be the property that, for every non empty family $\Cal D$ of polynomials, every special $\Cal D$-polynomial
$\Gamma$ of length $r \leq r_0$ is a sum of at most $C r_0 2^{r_0}$ gbs~$\zeta^{r}(\Cal D)$-polynomials, where $C$ is the constant
introduced before Lemma \ref{oneStep}.

Let $\Cal R =(R_j, j= 1, ..., d)$ be a family of $d$ polynomials. Let $\Gamma$ be an $\Cal R$-polynomial $\Gamma$ of length $r(\Gamma) = r_0+1$.
By applying the map $\Gamma \to \check \Gamma$ to $\Gamma$ which preserves the number of terms, we can assume that the $R_i$'s
are prime and distinct.

The property $H(2)$ is satisfied (the null polynomial is the only reduced special $\Cal R$-polynomial of length $\leq 2$), if the $R_i$'s
are pairwise relatively prime. Let us show that $H(r_0)$ implies $H(r_0 +1)$.

We use the fact that, if $\Gamma$ is a $p$th power of a special $\Cal R$-polynomial which is a sum of at most $t(r, \zeta^r(\Cal R))$
gbs~$\Cal S$-polynomials, then $\Gamma$ has the same property (with the same $t(r, \zeta^r(\Cal R))$) since the $p$th power of a sum is
the sum of the $p$th power of its terms.

Therefore, we can write $\Gamma = \x^\s \, \Lambda^{p^t}$, for some $\s \in \Z^d$ and some special $\Cal R$-polynomial $\Lambda$
(with the same number of terms: $r(\Lambda) = r(\Gamma)$) containing at least two non zero homogeneous components,
$\Lambda_{\be_0}, \Lambda_{\be_1}$. Multiplying by a monomial, one can assume $\be_0 = \0$. Let $\be= u(\be_1)$.

We apply Lemma \ref{oneStep} to $\Lambda$. With the previous notations, $\widetilde \Pi_\be(\Lambda)$ differs from
$\widetilde A_{\be}(\Lambda) + \widetilde B_\be(\Lambda)$ by at most $3 C r(\Gamma)$ gbs~$\zeta(\Cal R)$-polynomials.

$\widetilde A_{\be}(\Lambda)$ and  $\widetilde B_\be(\Lambda)$ are special
$\zeta(\Cal R)$-polynomials with a number of terms $\leq r_0$. Therefore, by the induction hypothesis (applied with $\Cal D = \zeta(\Cal R)$),
they are sum of at most $r_0 \, 2^{r_0}$ gbs~$\zeta^{r_0}(\zeta(\Cal R))$-polynomials.

For $r_0 \geq 3$, since $\zeta^{r_0}(\zeta(\Cal R)) = \zeta^{r_0+1}(\Cal R)$, $\widetilde \Pi$ is a sum of at most
$2 C r_0 \, 2^{r_0} + 3 C(r_0 + 1)\leq C (r_0+1) \, 2^{r_0+1}$ gbs~$\zeta^{r_0+1}(\Cal R)$-polynomials.

Using (\ref{decGamm1}), after multiplication of $\widetilde \Pi_\be(\Lambda)$ by $(\prod_{\rho \in \Upsilon(D_{\zeta^{r_0}(\Cal R), \beta, 1})}
x_\rho^{\theta_\rho(D_{\zeta^{r_0}(\Cal R), \beta, 1})})^{-1}$, the inverse of $\Psi(D_{r,\beta,1})$, to obtain $\Lambda$, this shows that $H(r_0+1)$
is true (a product of distinct prime polynomials is not a $p$-th power, hence its derivative is not zero, cf. Lemma \ref{nonZeroDer}).

The previous computation suffices to give an effective bound for the number of generalized basic special $\Cal E$-polynomials
in the decomposition of a polynomial $\Gamma$ of a given length $r(\Gamma)$. The following more precise estimation gives a polynomial bound.

First we take the $\be_0$-homogeneous component of $\Lambda$ which contains the biggest number of terms. Let $r \lambda_r$ be this number.
Altogether, the other components contain $r(1-  \lambda_r)$ terms.
Then we take the $\be_1$-homogeneous component which contains the second biggest number of terms (denoted by $r\mu_r$).

Let $c := p^{-d}$. As there are at most $p^d$ nonempty homogeneous components, we have $r\lambda_r \leq r-1$
(hence $1- \lambda_r \geq r^{-1}$) and $\lambda_r \geq c, \ \mu_r \geq c (1-  \lambda_r)$.

If $\theta > 2$ is such that $(1- c)^{\theta -1} \leq {c \over 2}$, then
$(1-  \lambda_r)^\theta + (1-  \mu_r)^\theta + {c \over 2} \, r^{1 - \theta} \leq 1$, since
\begin{eqnarray*}
&(1-  \lambda_r)^\theta + (1-  \mu_r)^\theta +  {c \over 2} r^{1 - \theta}\leq
{c \over 2} (1-  \lambda_r) + (1-  \mu_r)^\theta +  {c \over 2} r^{1 - \theta} \\\
&\leq {c \over 2} (1-  \lambda_r) + (1- \mu_r) + {c \over 2} r^{1 - \theta} \leq 1 + (1-  \lambda_r) (-{c \over 2}) +  {c \over 2} r^{1 - \theta}
\leq 1 - [{c \over 2} \, r^{-1} - {c \over 2} r^{1 - \theta}] \leq 1.
\end{eqnarray*}
For this choice of $\theta$ and $K = 6/c$, we have $K (r(1- \lambda_r))^\theta + K (r(1-  \mu_r))^\theta
+ 3 r \leq  K r^\theta$. Therefore this shows, by induction, that the number of needed gbs~$\zeta^r(\Cal R)$-polynomials
for the decomposition of $\Gamma$ is $\leq K r^\theta$. \eop

\vskip 2mm
\subsection{\bf Counting special $\Cal R$-polynomials}\label{counting}

\

We need an auxiliary lemma.
\begin{lem} \label{powp} Let $h$ be an integer, $F$ a finite set of non zero integers and $p$ an integer $> 1$.
For $h \geq 1$, let $W_h \subset \Z$ be the set of integers which can be written as a sum
$L = \sum_{i =1}^h v_i \, p^{t_i}$, $t_i \in \N$, $v_i \in F$.
There is a constant $K$ depending on $F, h$ such that, for all $N \geq 1$, the cardinal of the set $D_h \cap [-N, N]$
is less than $K\, (\log N)^h$.
\end{lem}
\proof Taking an element $L \not = 0$ in $W_h \cap [-N, N]$,
we can write $L = \sum_{i =1}^{h_1} v_i \, p^{t_i}$, $t_i \in \Z^+$, where we can assume that the set $\{t_j, j=1, ..., h_1\}$
is written in increasing order and $h_1 \leq h$ is such that $\sum_{i =k}^{h_1} v_i \, p^{t_i} \not = 0$, for all $1 \leq k \leq h_1$.

We have  $|L| = p^{t_1} \, |v_1 + \sum_{j=2}^{h_1} v_j \, p^{t_j-t_1}|$; hence $p^{t_1} \leq |L| \leq N$; therefore: $t_1 \leq \log N/\log p$.

Since $\sum_{j=2}^{h_1} v_j \, p^{t_j-t_1} \not = 0$, we have $1 \leq |\sum_{j=2}^{h_1} v_j \, p^{t_j-t_1}| \leq |L| + M$,
where $M$ denotes the maximum of $|u|$, for $u \in F$; hence: $p^{t_2-t_1} \, |v_2 + \sum_{j=3}^{h_1}  v_j \, p^{t_j-t_2}| \leq |L|+M$,
which implies
$$t_2 = t_2-t_1 + t_1 \leq \log (|L|+M) /\log p + \log |L| /\log p \leq \log (N+M) /\log p + \log N /\log p.$$
By iteration, we obtain $h_1 \leq h$ and a constant $C_h$ depending only on $h$ such that
$$L = \sum_{i =1}^{h_1} v_i \, p^{t_i} \text{ and } t_1 \leq t_2 \leq ... \leq t_{h_1} \leq C_h \log N /\log p.$$
Therefore $L$ can take at most $|2F|^h (C_h \log N)^h = K \, (\log N)^h$ different values. \eop

In the statement of the next theorem, $\Cal R$ is a family of polynomials $(R_1, ..., R_d)$ and $t(r)$ is the constant $t(r) = t(r, \Cal R)$
introduced in Theorem \ref{sumBasic0},
\begin{thm} \label{numbMinim} The number $\theta(D, r)$ of reduced special $\Cal R$-polynomials $\Gamma$ with $r$ terms,
supported in a domain $D$, satisfies for a constant $\gamma(r)$
\begin{eqnarray}
\theta(D, r) = O(|D|^{r/3} \, (\log \diam D)^{\gamma(r)}). \label{nbgons}
\end{eqnarray}
\end{thm}
\proof Let $\Gamma(\x) = \sum_{\a \in S(\Gamma)} c(\a) \, \x^\a$ be a reduced special $\Cal R$-polynomials with $r$ terms
such that $S(\Gamma) \subset D$.

By Theorem \ref{sumBasic0}, there are a finite family of polynomials $\Cal E$ and  $t = t(r, \Cal R)$ such that $\Gamma = \sum_{j=1}^{t} \Delta_j$,
where each $\Delta_j$ is a gbs~$\Cal E$-polynomials,
$$\Delta_j(\x) = \sum_{\b \in S(\Delta_j)} \, d(j,\b) \, \x^{\b}.$$
In the above formula, we have $\b \in \Z^{d'}$ for some $d^{\,'} \geq d$. If  $d^{\,'} > d$ we embed $\Z^d$ into $\Z^{d'}$
by completing by 0 the missing coordinated. We can view the elements $\a$ of the support $S(\Gamma)$ of $\Gamma$ as points in $\Z^{d'}$,
with the last $d^{\,'} - d$ coordinates equal to 0. The decomposition of $\Gamma$ reads more explicitly:
$$\Gamma(\x) = \sum_{\a \in S(\Gamma)} \, c(\a) \, \x^\a = \sum_j \, (\sum_{\b \in S(\Delta_j)} \, d(j,\b) \, \x^{\b})
=\sum_{\g \in \Z^{d'}} \, [\sum_j \sum_{\b \in S(\Delta_j): \, \b = \g} \, d(j,\b)] \, \x^\g.$$

Putting $u(\g) = \sum_j \sum_{\b \in S(\Delta_j) \, : \ \b = \g} \, d(j,\b)$, for $\g \in \Z^{d'}$, the formula reads in reduced form:
$$\Gamma(\x)= \sum_{\g \in \Z^{d'} \, : \, u(\g) \not = 0 } \, u(\g) \, \x^\g.$$
With the above embedding of $\Z^d$ into $\Z^{d'}$, we get
$$u(\g) = 0 \text{ for } \g \not \in S(\Gamma), \ u(\g) = c(\g) \text{ for } \g \in S(\Gamma).$$

Let us denote by $\Phi$ the family of the $\Delta_j$'s and write $d(\Delta, \b)$ instead of $d(j,\b)$ for the coefficients
of $\Delta = \Delta_j \in \Phi$.

Using an idea of \cite{ArBeBe08}, we put a graph structure on $\Phi$ by saying that there is an arrow between $\Delta$ and $\Delta'$ if
$S(\Delta) \cap S(\Delta') \not = \emptyset$. For this graph structure, $\Phi$ decomposes in connected components denoted by $\Phi_k$.

Let $S_k:= \bigcup_{\Delta \in \Phi_k} S(\Delta)$. Since $S(\Delta)$ and $S(\Delta')$ are disjoint for $\Delta, \Delta'$ in different components,
the sets $S_k$ are disjoint.

It follows that the above definition of $u$  can be written
$$u(\g) = \sum_{\Delta \in \Phi_k} \, \sum_{\b \in S(\Delta) \, : \ \b = \g}  \, d(\Delta,\b), \text{ for } \g \in S_k.$$

We have $\Gamma(\x) = \sum_k \Gamma_k(\x)$ with
\begin{eqnarray*}
\Gamma_k(\x) &&= \sum_{\g \in S_k} \, u(\g) \, \x^\g = \sum_{\g \in S_k} \,
[\sum_{\Delta \in \Phi_k} \, \sum_{\b \in S(\Delta) \, : \ \b = \g}  \, d(\Delta,\b)] \, \x^\g \\
&&=\sum_{\Delta \in \Phi_k} \, [\sum_{\b \in S(\Delta)}  \, d(\Delta,\b)] \, \x^\b] = \sum_{\Delta \in \Phi_k} \, \Delta(\x).
\end{eqnarray*}
Therefore the supports $S(\Gamma_k)$ are pairwise disjoint and each $\Gamma_k$ is a special $\Cal E$-polynomial
(actually, once reduced, a special $\Cal R$-polynomial).

We say that a reduced special $\Cal R$-polynomial $\Gamma = \sum_{\a \in S(\Gamma)} c(\a) \, \x^\a$ is  $\Cal R$-{\it minimal} if,
for every set $S_1$ strictly contained in $S(\Gamma)$, the polynomial $\sum_{\a \in S_1} c(\a) \, \x^\a$ is not a special $\Cal R$-polynomial.

Let us assume first that the polynomial $\Gamma$ is $\Cal R$-minimal. The disjointness of the supports $S(\Gamma_k)$  implies that $\Phi$ is
a connected graph.

The support of the gbs~polynomials are sites of the form $\b + p^k \v_t, t \in J$, where $J$ is a finite set of indices corresponding
to the collection of all basic special polynomials. Suppose that $\Delta, \Delta'$ are two gbs~polynomials with a common site in their support.
This site reads
$\b + p^k \v_t = \b' + p^{k'} \v_{t'}$, since it belongs to $\Delta$ and $\Delta'$. Therefore, $\b' - \b = p^{k} \v_t - p^{k'} \v_{t'}$.

If $\c_1$ and $\c_2$ belong respectively to $\Delta$ and $\Delta'$, then we have:
$\c_1  = \b + p^{k_1} \v_{t_1}$, $\c_2 = \b' + p^{k_2} \v_{t_2}$;
hence: $\c_2 - \c_1 = \b' + p^{k_2} \v_{t_2} - (\b + p^{k_1} \v_{t_1}) = p^{k} \v_t - p^{k'} \v_{t'} + p^{k_2} \v_{t_2} - p^{k_1} \v_{t_1}$.

It follows that, if $\c'$ belongs to a connected chain (starting at $\c$) of gbs~polynomials $\Delta_j$
(i.e., two consecutive $\Delta, \Delta'$ in the chain have a common site in their support), the difference $\c' - \c $ has the form:
\begin{eqnarray}
\c' - \c = \sum_i (p^{k_i} \v_{t_i} - p^{k_i'} \v_{t_i'}). \label{connectChn}
\end{eqnarray}
There are $|D|$ choices for $\c$ multiplied by $p-1$ (the cardinal of $\F_p \setminus \{0\}$). We obtain all minimal special $\Cal R$-polynomials
starting from $\c$ by constructing all possible connected chains of gbs~$\Cal E$-polynomials.

Since, in view of (\ref{connectChn}), $S(\Gamma) \subset \prod_{i=1}^d (W_{2t(r)} + c_i)$, using Lemma \ref{powp} for each coordinate, we obtain
that the number of choices is at most, for a given starting point $\c$, $[K \, (\log \diam(D))^{t_1(r)}]^d$, where $t_1(r)$ is a constant.

This implies that the number $\theta(D, r)$ of minimal special $\Cal R$-polynomials $\Gamma$ with $s$ terms, $s \leq r$,
supported in a domain $D$, satisfies the bound \begin{eqnarray}
\theta(D, r) = O(|D| \, (\log \diam D)^{d \, t_1(r)}). \label{nbMingons}
\end{eqnarray}
If $\Gamma$ is not $\Cal R$-minimal, then there is $S_1$ strictly contained in $S(\Gamma)$ such that $\sum_{\a \in S_1} c(\a) \, \x^\a$
is a special $\Cal R$-polynomial. Since $\sum_{\a \in S(\Gamma) \setminus S_1} c(\a) \, \x^\a$ is also a special $\Cal R$-polynomial,
by iteration of this decomposition, any special $\Cal D$-polynomial decomposes as a sum of minimal ones with disjoint supports.
As the length of a minimal polynomial is at least 3, (\ref{nbgons}) follows from (\ref{nbMingons}).
\eop

\section{\bf Application to limit theorems} \label{almMixSect}

\subsection{\bf Preliminaries: variance, cumulants}\label{prelimAppl}

\

We need some general facts about variance, summation sequences, cumulants. (See \cite{CohCo15} for more details.)
Recall that, if ${\Cal S} = (T^\el, \el \in \Z^d)$ is an abelian group isomorphic to $\Z^d$ of unitary
operators on a Hilbert space $\Cal H$, for every $f \in \cH$ there is a positive finite measure
$\nu_f$ on $\T^d$, the spectral measure of $f$, with Fourier
coefficients $\widehat \nu_f({\el}) = \langle T^{\el} f, f \rangle$,
$\el \in \Z^d$. When $\nu_f$ is absolutely continuous, its density
is denoted by $\varphi_f$.

We assume that $\Cal S$  has the Lebesgue spectrum property for its action on ${\Cal H}$, i.e.,
there exists a closed subspace ${\Cal K}_0$ such that $\{T^\el {\Cal K}_0, \, \el \in \Z^d\}$ is a family of pairwise orthogonal
subspaces spanning a dense subspace in ${\Cal H}$. If $(\psi_j)_{j \in \Cal J}$ is an orthonormal basis of ${\Cal K}_0$, $\{T^\el
\psi_j, j \in \Cal J, \, \el \in \Z^d \}$ is an orthonormal basis of ${\Cal H}$.
For every $f \in \cH$, $\nu_f$ has a density $\varphi_f$ in $L^1(d\t)$.

{\it Summation sequence}

{\it Definitions}: We call {\it summation sequence} any sequence $(w_n)_{n \geq 1}$ of functions from $\Z^d$ to $\R^+$
with $0 < \sum_{\el \in \Z^d} w_n(\el) < +\infty, \ \forall n \geq 1$. Given ${\Cal S} = \{T^\el, \el \in \Z^d\}$ and $f \in {\Cal H}$, the
associated sums are $\sum_{\el \in \Z^d} w_n(\el) \, T^\el f$.

We say that $(w_n)$ is $\zeta$-{\it regular}, if $\zeta$ is a probability measure on $\T^d$ and the sequence of nonnegative kernel $\tilde w_n$ defined by
\begin{eqnarray}
\tilde w_n(\t) = {|\sum_{\el \in \Z^d} w_n(\el) \ e^{2\pi i \langle
\el, \t \rangle}|^2 \over \sum_{\el \in \Z^d} |w_n(\el)|^2}, \ \t \in \T^d,\label{chch}
\end{eqnarray} weakly converges to $\zeta$ when $n$ tends to infinity. This is equivalent to
$$\widehat \zeta (\p) = \lim_{n \to \infty} \int \tilde w_n(\t) \, e^{-2\pi i \langle \p, \t \rangle} \ d\t, \forall \p \in \Z^d.$$
When the spectral density is continuous, $\varphi_f \to (\zeta(\varphi_f))^\frac12$ satisfies the triangular inequality.

{\it Variance for summation sequences}

If $(w_n)$ is a $\zeta$-regular summation sequence and $f$ in $\Cal H$ with a continuous spectral density
$\varphi_f$. By the spectral theorem, we have for $\underline \theta = (\theta_1, ..., \theta_d) \in \T^d$:
\begin{eqnarray}
(\sum_\el \, w_n^2(\el))^{-1} \|\sum_{\el \in \Z^d} \, w_n(\el) \,
e^{2\pi i \langle \el, \underline \theta \rangle} \, T^\el f\|_2^2 = (\tilde
w_n * \varphi_f)(\underline \theta) \underset{n \to \infty} {\longrightarrow}
(\zeta * \varphi_f)(\underline \theta). \label{spectrThm}
\end{eqnarray}
For example, if $(D_n)$ is a F\o{}lner sequence of sets in $\Z^d$, then $w_n(\el) = 1_{D_n}(\el)$, $\zeta = \delta_\0$ and the
usual asymptotic variance $\sigma^2(f)$ is $\varphi_f(\0)$.

\goodbreak
{\it Moments, cumulants and the CLT}

Let us recall now some general results on mixing of order $r$, moments and cumulants (see \cite{Leo60a}).
In what follows, we assume the random variables to be uniformly bounded.

Let $(X_1, ... , X_r)$ be a random vector. For any subset $I = \{i_1, ..., i_p\} \subset J_r:= \{1, ..., r\}$, we put
$m(I) = m(i_1, ..., i_p):= \E(X_{i_1} \cdots X_{i_p})$.
Cumulants are computed from moments by
\begin{eqnarray}
C(X_1, ... , X_r) = \sum_{\pi \in \Cal P} (-1)^{p(\pi)-1} (p(\pi) - 1)! \, \,  m(I_1) \cdots m(I_{p(\pi)}), \label{cumForm0}
\end{eqnarray}
where $\pi = \{I_1, I_2, ..., I_{p(\pi)}\}$ runs through the set $\Cal P$ of partitions of $J_r = \{1, ..., r\}$ into nonempty subsets
and $p(\pi)$ is the number of elements of $\pi$.

Putting $s(I) := C(X_{i_1}, ... , X_{i_p})$ for $I = \{i_1, ..., i_p\}$, we have
\begin{eqnarray}
\E(X_{1} \cdots X_{r}) = \sum_{\pi \in \Cal P} s(I_1) \cdots s(I_{p(\pi)}).
\label{cumFormu2}
\end{eqnarray}
For a single random variable $Y$, we define $C^{(r)}(Y) := C(Y, ..., Y)$, where $(Y, ..., Y)$ is the vector with $r$ components
equal to $Y$. If $Y$ is centered, $C^{(2)}(Y)$ coincides with $\|Y\|_2^2$.

Let be given a random field of real random variables $(X_\k)_{\k \in \Z^d}$ and a summable weight $w$ from $\Z^d$ to $\R$. For $Y
:= \sum_{\el \in \Z^d} \, w(\el) \, X_\el$, using the multilinearity of the cumulants, we obtain:
\begin{eqnarray}
C^{(r)}(Y) = \sum_{(\el_1, ..., \el_r) \, \in (\Z^d)^r} \,  w(\el_1) \cdots w(\el_r)  \, C(X_{\el_1}, \cdots , X_{\el_r}) \,. \label{cumLin0}
\end{eqnarray}
\begin{lem} The number $\gamma(p, r)$ of partitions of $J_r$ into $p \le r$ nonempty subsets satisfies
\begin{eqnarray}
\sum_{p=1}^r (-1)^{p-1} \, (p - 1)! \ \gamma(p, r) = 0. \label{sumzero}
\end{eqnarray}
\end{lem}
\proof  (\ref{sumzero}) follows by induction from the following formula:
$\gamma(p, r) = \gamma(p-1, r-1) + \, p \, \gamma(p, r-1), \ p = 1, ..., r, \, r \geq 1$.
\eop

\begin{thm} \label{Leonv} (cf. \cite{Leo60b}, Theorem 7)
Let $(X_\k)_{\k\in \Z^d}$ be a random process and $(w_n)_{n \geq 1}$
a summation sequence on $\Z^d$. Let $Y^{(n)} = \sum_\el \, w_n(\el) \, X_\el, \ n \geq 1$. If $\|Y^{(n)}\|_2 \not = 0$ and
\begin{eqnarray}
\sum_{(\el_1, ..., \el_r) \, \in (\Z^d)^r} \, w_n(\el_1) ... w_n(\el_r) \, C(X_{\el_1}, ... , X_{\el_r})
= o(\|Y^{(n)}\|_2^r), \forall r \geq 3, \label{smallCumul}
\end{eqnarray}
then ${Y^{(n)} \over \|Y^{(n)}\|_2}$ tends in distribution to $\Cal N(0, 1)$
when $n$ tends to $\infty$.
\end{thm}
\proof \ Let $\beta_n := \|Y^{(n)}\|_2 = \|\sum_\el \, w_n(\el) \, X_\el\|_2$ and $Z^{(n)} = \beta_n^{-1} Y^{(n)}$.
In view of (\ref{cumLin0}), we have $C^{(r)}(Z^{(n)}) = \beta_n^{-r}
\sum_{(\el_1, ..., \el_r) \, \in (\Z^d)^r} \, w(\el_1) ... w(\el_r) \, C(X_{\el_1}, ... , X_{\el_r})$, hence by (\ref{smallCumul}):
\begin{eqnarray}
\lim_n C^{(2)}(Z^{(n)}) = 1, \ \lim_n C^{(r)}(Z^{(n)}) = 0, \forall r \geq 3. \label{limCum}
\end{eqnarray}
Using the formula linking moments and cumulants, the theorem follows from the result of \cite{FreSho31}
applied to $(Z^{(n)})_{n \geq 1}$. \eop

\vskip 3mm
{\it Algebraic framework}

Coming back to the framework of a compact abelian group $G$, we consider a totally ergodic $\N^d$-action $\el \to T^\el$
by algebraic commuting endomorphisms on $G$, or its invertible $\Z^d$-extension, with the Lebesgue spectrum property.

Below a function $f$ on $G$ will be called a "regular function" if $f$ belongs to the space $AC_0(G)$, i.e., has an absolutely convergent
Fourier series. Recall that, if $f$ is regular, its spectral density $\varphi_f$ is continuous on $\T^d$
and for every $\varepsilon > 0$ there is a trigonometric polynomial $P$ defined on $G$ such that $\|\varphi_{f - P}\|_\infty \leq \varepsilon$.

The proof of the CLT given in \cite{Leo60b} for a single ergodic endomorphism of a compact abelian group $G$ is based on
the computation of the moments of the ergodic sums of trigonometric polynomials and uses mixing of all orders.
As mentioned in Section \ref{actFp}, for $\Z^d$-actions by automorphisms on $G$, mixing of all orders is satisfied when $G$ is connected, but may fail
for non connected groups like shift-invariant subgroups of $\F_p^{Z^d}$. Nevertheless, when the non-mixing configurations are sparse enough,
the moment method can be applied.

{\it Non-mixing $r$-tuples}

Let $f = \sum_{j \in J} c_j \chi_j$ be a trigonometric polynomial and $\Phi = (\chi_j, j \in J)$.
We defined the set of ``{\it non-mixing}" $r$-tuples for $\Phi = (\chi_j, j \in J)$ by
\begin{eqnarray}
&&{\cal N}(\Phi, r) := \{(\a_1, ..., \a_r): \, \exists \chi_{j_1}, ... \chi_{j_r} \in \Phi:
C(T^{\a_1} \chi_{j_1}, ..., T^{\a_r} \chi_{j_r}) \not = 0\}. \label{non-mixingGon}
\end{eqnarray}
In view of (\ref{r-moment}) (appendix) and (\ref{cumForm0}), if $(\a_1, ..., \a_r) \in {\cal N}(\Phi, r)$,
we have $T^{\a_1} \chi_{j_1} ...  T^{\a_r} \chi_{j_r} = \chi_0$, for some $(\chi_{j_1}, ..., \chi_{j_r}) \in \Phi$.
We will use the results of the subsection \ref{counting} to show that the sets ${\cal N}(\Phi, r)$ are small in some sense.

\vskip 3mm
\subsection{\bf Counting non zero cumulants}

\

Now we consider the action by endomorphisms discussed in the first section.
For $d \geq 2$, $R_1, R_2, ..., R_d$ are $d$ polynomials of degree $\geq 1$ in $x$ over $\F_p$, fixed once for all.
Recall that for $\a_j \in \Z^d$, the action of $T^{\a_j}$ on a character $\chi_{Q_j}$ associated to a polynomial $Q_j$
is the multiplication of $Q_j$ by $\uR^{\a_j}= \prod_{i=1}^d R_i^{a_{j,i}}$.

For $\widetilde Q = (Q_1, ..., Q_r)$, the corresponding cumulant is
$C_{\widetilde Q}(A) = C(T^{\a_1} \chi_{Q_1}, ..., T^{\a_r} \chi_{Q_r})$.
Let $\chi_1, ..., \chi_r$ be characters on $\F_p^{\Z}$. They correspond to a set of polynomials in one variable $\widetilde Q = \{Q_1, ..., Q_r\}$.
For an $r$-tuple $A = (\a_1, ..., \a_r) \in (\Z^d)^r$ the relation $T^{\a_1} \chi_1  ... T^{\a_r} \chi_r = \chi_0$ is equivalent to the relation
\begin{eqnarray}
\sum_{j= 1}^r Q_j \, \prod_{i=1}^d R_i^{a_{j,i}} = 0. \label{equaCharcB}
\end{eqnarray}

In the present framework, the formula for cumulants is used for the random variables $X_j= T^{\a_j} \chi_j$, where the characters $\chi_j = \chi_{Q_j}$
are associated by (\ref{characp}) to {\it non zero given fixed polynomials (over $\F_p$) $Q_j$, $i=1, ...,r$}.

For a domain $D \subset\Z^d$, $D^r$ denotes the set of $r$-tuples $\uA$ of elements of $D$.

Let $Q := \sum_i Q_i \uR^{\a_i}$. The moments read as the integral (actually a finite discrete sum)
$$\int e^{{2\pi \over p}i \, \sum_{k \in S(Q)} c(Q, k) \, \zeta_k} \, d\zeta =\prod_{k \in S(Q)} \frac1p \, \sum_{j=0}^{p-1}e^{{2\pi \over p}i \,
c(Q, k) j}.$$
They are equal to 1 if $\sum_i Q_i \uR^{\a_i} = 0$ and to 0 else (mod $p$).
\begin{proposition} \label{bndzerocuml} For each $r \geq 3$, there are constants $\gamma, K$ (dependent on $\widetilde Q$)
such that
\begin{eqnarray}
\# \{\uA \in D^r : C_{\widetilde Q}(\uA) \not = 0\} \leq K \, |D|^{{r\over 2} - \frac12} \, (\log \diam D)^{\gamma}. \label{majnbrCm}
\end{eqnarray}
\end{proposition}
\proof If $C(T^{\a_1} \chi_1, ..., T^{\a_r} \chi_r) \not = 0$, by (\ref{cumForm0}) there exists a partition $\pi = \{I_1, ..., I_p\}$
of $J= \{1, ..., r\}$ such that $\sum_{j \in I_k} \, Q_j \, \uR^{\a_j} = 0, \ k=1, ..., p.$
This implies $\sum_{j \in J} \, Q_j \, \uR^{\a_j} = 0$.
The polynomial $\Lambda(x, \x) = \sum_{j= 1}^r Q_j(x) \, \prod_{i=1}^d x_i^{a_{j,i}}$ satisfies (\ref{equaCharcB}) when $R_i$ is substituted to $x_i$.

Let $\Upsilon(\widetilde Q)$ be the set of prime factors of the polynomials $Q_j$ in $\widetilde Q$.
In $\Upsilon(\widetilde Q)$ it may exist prime factors belonging to $\Cal R$ and possibly new prime factors denoted by
$R_i$, $i =d+1, ..., \delta$. We enlarge the set $\cal R$ to $\widetilde {\cal R} = \cal R \bigcup \Upsilon(\widetilde Q)$ by adding to $\cal R$
the prime factors of the $Q_j$'s, i.e., we consider the set of prime polynomials $\widetilde {\cal R} = \{R_1, ..., R_d, R_{d+1}, ..., R_{d+\delta}\}$.

The factorization of $Q_j$ in prime monic polynomials (with $d(Q_j) \in \F_p$) is
$$Q_j(x) = d(Q_j) \prod_{\rho \in \Upsilon(\widetilde Q)} \, \rho^{g_{j,\rho}}
= d(Q_j) \, \prod_{i=1}^d \, R_i^{g_{j,i}} \, \prod_{i=d+1}^{d+\delta} \, R_i^{g_{j,i}}.$$
Some of the $g_{j,i}$ may be zero. Equation (\ref{equaCharcB}) reads
$$\sum_{j= 1}^r  \, d(Q_j) \, \prod_{i=1}^d \, R_i^{g_{j,i}} \, \prod_{i=d+1}^{d+\delta} \, R_i^{g_{j,i}} \prod_{i=1}^d R_i^{a_{j,i}} = 0.$$
Putting $a_{j,i} = 0$ for $i=d+1, ..., d+\delta$, the new $r$-tuple
$B = (\b_1, ..., \b_r)$ in $(\Z^{d'})^r$ is given by $b_{j,i}= a_{j,i} + g_{j,i}$, $i=1, ..., d+\delta$.
We get a polynomial with $d' \geq d$ variables, $\sum_{j= 1}^r\, c(\b_j) \, \prod_i^{d+\delta} x_i^{b_{j,i}}$
which is a (not necessarily reduced) special $\widetilde {\cal R}$-polynomial. We have
$$\sum_{j \in J} \, Q_j \, \uR^{\a_j} = \sum_{\b} \, c(\b) \, \widetilde \uR^{\b}, \text{ with } c(\b) = \sum_{j: \,\a_j+\g_j = \b} d(Q_j).$$
The $r$-tuple $\uA$ can be viewed as a collection of $r$ vectors in $\Z^d$ which is divided into the two following subsets:
$\uA_0:= \{\a_j: \, c(\a_j + g_j) = 0\}$, $\uA_1:= \{\a_j: \, c(\a_j + g_j) \not = 0\}$.
The terms corresponding to $\a_j \in \uA_0$ disappear. The sum $\Gamma(\x) = \sum_{\b: \, c(\b) \not = 0} \, c(\b) \, \x^\b$ is reduced.

Once the sets $\uA_0, \uA_1$ are chosen, $\uA$ is determined up to a permutation which introduces a bounded factor in the counting
of the configurations $\uA_0$.

In what follows, $K$ will a generic constant which may change from an inequality to another. $\widetilde D$ is the domain obtained from $D$
when the $a$'s are replaced by the $b$'s. Its cardinal and its diameter are less than a constant times the
cardinal and the diameter of $D$.

Let us say that $\a_j$ is equivalent to $\a_{j'}$ if $\a_{j} + \g_{j'} = \a_{j'} + \g_{j'}$.
All elements in the same equivalence class are at bounded distance from each other (their mutual distance is bounded by $\max_{j,j'} \|g_j - g_{j'}\|$.
Once an element is chosen in a class, there is an uniformly bounded number of choices for the other elements.
The classes of $\a_j$'s such that $c(\a_j +\g_j) = 0$ have at least two elements.

Let $t \in [0, r]$ be the number of elements in $\uA_0$. The number of choices for the elements of $\uA$ belonging to $\uA_0$ is at most
$K\, |D|^{t/2}$. The polynomial $\Gamma$ is reduced and has less than $r-t$ terms. By Theorem \ref{numbMinim}
the number of choices of such polynomials is less than $K \, |\tilde D|^{(r-t)/3} \, (\log \diam \tilde D)^{\gamma(r)}$.

Therefore the total number of choices for $\uA$ is at most $K \, |\tilde D|^{t/2 + (r-t)/3} \, (\log \diam \tilde D)^{\gamma(r)}
= K \, |\tilde D|^{r/3 + t/6} \, (\log \diam \tilde D)^{\gamma(r)}$.

If $\uA_1$ is not empty, then we have $r-t \geq 3$, since a reduced special $\Cal R$-polynomial has at least 3 terms, and the above upper bound is
less than $K \, |\tilde D|^{r/2 - 1/2} \, (\log \diam \tilde D)^{\gamma(r)}$.

If $\uA_1$ is empty, then $t = r$. If each class is composed only of pairs of 2 elements, then $r=2r'$ is even
and the computation of the cumulant corresponds exactly (for $r'$ instead of $r$) to the case where all moments are equal to 1.
By (\ref{sumzero}) the cumulant is 0. It shows that this case does not appear in the computation for (\ref{majnbrCm}).
Therefore there is a class containing at least 3 elements and we have a bound by
$K \, |\tilde D|^{(r-3)/2 + 1} \, (\log \diam \tilde D)^{\theta(r)}) \leq  K \, |\tilde D|^{r/2 - 1/2} \, (\log \diam \tilde D)^{\gamma(r)})$.
which is less than  $K \, |D|^{r/2 - 1/2} \, (\log \diam D)^{\gamma(r)})$, for a new constant $K$.
\eop

\vskip 2mm
{\it Example:} For $r=4$, the cumulants are given by
\begin{eqnarray*}
&\int T^{\a_1} \chi_1 \, T^{\a_2} \chi_2 \, T^{\a_3} \chi_3 \, T^{\a_4} \chi_4 \,
- [\int T^{\a_1} \chi_1 \, T^{\a_2} \chi_2 \,
\int T^{\a_3} \chi_3 \, T^{\a_4} \chi_4 \,\\
&+ \int T^{\a_1} \chi_1 \, T^{\a_3} \chi_3 \,
 \int T^{\a_2} \chi_2 \, T^{\a_4} \chi_4 \, + \int T^{\a_1} \chi_1 \, T^{\a_4} \chi_4 \,
 \int T^{\a_2} \chi_2 \, T^{\a_3} \chi_3].
\end{eqnarray*}

The characters are given by polynomials $Q_i$. The integrals and their products take the value 0 or 1. Each time an integral is 1,
we have relations of the form $\sum_{i \in I} Q_i \, R^{\a_i}  = 0$. There are 3 cases:

a) $Q_1 \, R^{\a_1} + Q_2 \, R^{\a_2} + Q_3 \, R^{\a_3} + Q_4 \, R^{\a_4} = 0$ (and no vanishing subsums), \\
b) $Q_1 \, R^{\a_1} + Q_2 \, R^{\a_2} = 0$ and $Q_3 \, R^{\a_3} + Q_4 \, R^{\a_4} = 0$, or the analogous relations obtained by permutation, \\
c) [$Q_1 \, R^{\a_1} + Q_2 \, R^{\a_2} = 0$, $Q_3 \, R^{\a_3} + Q_4 \, R^{\a_4} = 0$],
[$Q_1 \, R^{\a_1} + Q_3 \, R^{\a_3} = 0$ and $Q_2 \, R^{\a_2} + Q_4 \, R^{\a_4} = 0$],
or the analogous relations obtained by permutation.

In case c) we see that $\a_1$ and $\a_2$ are close together as well as $\a_3$ and $\a_4$ and $\a_2$ and $\a_4$.
It follows that the four elements $\a_1, \a_2, \a_3, \a_4$ are close together and there is only one degree of freedom
for the choice of $A$ if $A$ belongs to this type of 4-tuple.

If we are in case b), but not in case c), then we have the relations $Q_1 \, R^{\a_1} + Q_2 \, R^{\a_2} = 0$
and $Q_3 \, R^{\a_3} + Q_4 \, R^{\a_4} = 0$ (hence $Q_1 \, R^{\a_1} + Q_2 \, R^{\a_2} + Q_3 \, R^{\a_3} + Q_4 \, R^{\a_4} = 0$).
The cumulant reduces to $\int T^{\a_1} \chi_1 \, T^{\a_2} \chi_2 \, T^{\a_3} \chi_3 \, T^{\a_4} \chi_4 \,
- \int T^{\a_1} \chi_1 \, T^{\a_2} \chi_2 \, \int T^{\a_3} \chi_3 \, T^{\a_4} \chi_4 \,= 1 -1 = 0$.

If we are in case a), but not b) or c), the cumulant is $\int T^{\a_1} \chi_1 \, T^{\a_2} \chi_2 \, T^{\a_3} \chi_3 \, T^{\a_4} \chi_4 \, = 1$.
The relation $Q_1 \, R^{\a_1} + Q_2 \, R^{\a_2} + Q_3 \, R^{\a_3} + Q_4 \, R^{\a_4} = 0$ may be reducible, but we find at least 3 terms in the irreducible
relations of the decomposition. The number of 4-tuples $A = (\a_1, \a_2, \a_3, \a_4)$ belonging to types corresponding to case a) is
less than $O(|D| \, (\log \diam D)^\theta)$ for some constant $\theta$.

\vskip 3mm
\subsection{\bf Examples of limit theorems for some shift-invariant groups} \label{quenchSect}

\

If $(w_n)_{n \geq 1}$ is a summation sequence on $\Z^d$, for $f \in L^2(G)$, we put $\sigma_n(f) := \|\sum_\el w_n(\el) \, T^\el f \|_2$
and assume $\sigma_n^2(f) \not = 0$, for $n$ big enough. We suppose that $(w_n)$ is $\zeta$-regular.

We can suppose $\zeta(\varphi_f) > 0$, since otherwise the limiting distribution is $\delta_0$.
By $\zeta$-regularity we have $\sigma_n^2(f) \sim (\sum_{\el} \, w_n^2(\el)) \, \zeta(\varphi_{f})$ with $\zeta(\varphi_{f}) > 0$.
\begin{thm} \label{tclRegKer} Let $(w_n)_{n \geq 1}$ be a summation sequence on $\Z^d$ which is $\zeta$-regular
(cf. definition in Subsection \ref{prelimAppl}).
Let $f$ be a regular function with spectral density $\varphi_f$ such that  $\zeta(\varphi_f) > 0$. The condition
\begin{eqnarray}
&& \sum_{(\el_1, ..., \el_r) \in {\cal N}(\Phi,r)} \prod_{j= 1}^r w_n(\el_j) = o\bigl((\sum_{\el \in \Z^d}
w_n^2(\el))^\frac{r}2\bigr), \, \forall \text{ finite family } \Phi \text{ of characters }, \forall r \geq 3, \label{cumulCNab}
\end{eqnarray}
implies
\begin{eqnarray}
(\sum_{\el \in \Z^d} \, w_n^2(\el))^{-\frac12} \, \sum_{\el \in
\Z^d} w_n(\el) f(T^\el .) \overset{distr} {\underset{n \to \infty}
\longrightarrow } \Cal N(0, \zeta(\varphi_f)). \label{cvgce1}
\end{eqnarray}
\end{thm}
\proof
a) First let us take for $f$ a trigonometric polynomial.
Let us check (\ref{smallCumul}) of Theorem \ref{Leonv}, i.e., in view of (\ref{cumLin0}),
\begin{eqnarray}
|\sum_{(\el_1, ..., \el_r) \, \in \, (\Z^d)^r} C(T^{\el_1} f,..., T^{\el_r} f) \, w_n(\el
_1)...w_n(\el _r)| = o\bigl((\sum_{\el \in \Z^d} w_n^2(\el))^\frac{r}2\bigr), \, \forall r \geq 3. \label{condCumul1}
\end{eqnarray}
If the cumulant $C(T^{\el_1} f,..., T^{\el_r} f)$ is $\not = 0$, then $(\el _1, ..., \el _r)$ is a non-mixing $r$-tuple
for the set $\Phi$ of characters which appear in the expansion of $f$; hence (cf. notation (\ref{non-mixingGon})):
\begin{eqnarray*}
\sum_{(\el_1, ..., \el_r) \in (\Z^d)^r} \, C(T^{\el_1} f,..., T^{\el_r} f) \, \prod_{j= 1}^r w_n(\el_j) =
\sum_{(\el_1, ..., \el_r) \in {\cal N}(\Phi,r)} \,  C(T^{\el_1} f,..., T^{\el_r} f) \, \prod_{j = 1}^r w_n(\el_j).
\end{eqnarray*}
Since the cumulants are bounded, the sums in the previous formula are bounded by $C \sum_{(\el_1, ..., \el_r)
\in {\cal N}(\Phi,r)} \prod_{j= 1}^r w_n(\el_j)$. Therefore, in view of (\ref{cumulCNab}), the condition of Theorem \ref{Leonv} is satisfied.
This implies the CLT when $f$ is a trigonometric polynomial.

b) Now, for a regular function by the $\zeta$-regularity of $(w_n)$, we have:
$$(\sum_{\el \in \Z^d} \, w_n^2(\el))^{-1} \, \|\sum_{\el \in \Z^d} w_n(\el)
\, T^\el f \|_2^2 = \int_{\T^d} \, \tilde w_n \, \varphi_f \, dt
\underset{n \to \infty} \to \zeta(\varphi_f).$$
If $(\varepsilon_k)$ a sequence of positive numbers tending to 0, there is a sequence of trigonometric polynomials $(f_k)$ such that:
$\|\varphi_{f - f_k}\|_\infty \leq \varepsilon_k$. Let us consider the processes defined respectively by
\begin{eqnarray*}
U_n^{(k)} := (\sum_{\el \in \Z^d} \, w_n^2(\el))^{-\frac12} \,
\sum_{\el \in \Z^d} w_n(\el) \, f_k(T^\el .), \ U_n := (\sum_{\el \in
\Z^d} \, w_n^2(\el))^{-\frac12} \, \sum_{\el \in \Z^d} w_n(\el) \, f(T^\el .).
\end{eqnarray*}
We have $\zeta(\varphi_{f - f_k}) \to 0$. It follows that $\zeta(\varphi_{f_k}) \not = 0$ for $k$ big enough.

Since $\sigma_n^2(f_k) \sim (\sum_{\el} \, w_n^2(\el)) \, \zeta(\varphi_{f_k})$ with $\zeta(\varphi_{f_k}) > 0$,
it follows from the result in a) for the trigonometric polynomials $f_k$: ${U_n^{(k)} \overset{distr} {\underset{n \to \infty}
\longrightarrow} \Cal N(0,\zeta(\varphi_{f_k}))}$ for every fixed $k$. Moreover, since
\begin{eqnarray*}
\lim_n \int |U_n^{(k)} - U_n|_2^2 \ d\mu &=& \lim_n \int_{\T^d} \,
\tilde w_n \,\varphi_{f-f_k} \, d \t = \zeta(\varphi_{f - f_k}) \leq \varepsilon_k,
\end{eqnarray*}
we have $\limsup_n \mu[|U_n{(k)} - U_n| > \varepsilon] \leq \varepsilon^{-2} \limsup_n \int |U_n{(k)} - U_n|_2^2 \ d\mu \underset{k
\to \infty} \to 0$ for every $\varepsilon > 0$.

Therefore the condition $\lim_k \limsup_n \mu[|U_n{(k)} - U_n| > \varepsilon] = 0$, $\forall \varepsilon > 0$, is satisfied and the
conclusion $U_n \overset{distr} {\underset{n \to \infty} \longrightarrow } \Cal N(0, \zeta(\varphi_{f}))$ follows from
Theorem 3.2 in \cite{Bill99}. \eop

\vskip 3mm
{\bf Application to shift-invariant subgroups}

The limit theorems shown in \cite{CohCo15} hold in the present framework of shift-invariant subgroups.
We restrict the presentation to two examples.

Let us consider a family $(R_j, j \in J)$ of polynomials of degree $\geq 1$ and $\gamma_{j} = \gamma_{R_j}$ the corresponding endomorphisms
of $K = \F_p^{\Z}$. As in Section \ref{shiftInv}, taking the natural invertible extension, we extend them to automorphisms of the
shift-invariant subgroup $G_{\Cal J }$ of $G^{(d+1)}$ defined by the ideal $\Cal J = \Ker(h_{\cal R})$. The $(R_j)$'s are chosen to be algebraically
independent. Therefore we have a totally ergodic $\Z^d$-action $(T^\el, \el \in \Z^d)$ on
$G_{\Cal J}$, with $T^\el = T_1^{\ell_1}... T_d^{\ell_d}$ and $T_j$ the composition by the shift $\sigma_{j+1}$.

\goodbreak
{\bf Example 1: F\o{}lner sequence in $\N^d$}

\begin{thm} \label{tclFoln} Let $(D_n)_{n \geq 1}$ be a F\o{}lner sequence of sets in $\N^d$. If $f$ is a regular function,
we have $\sigma^2(f) = \lim_n \|\sum_{\el \in D_n} \, T^\el f\|_2^2/|D_n| = \varphi_f(0)$. If moreover $\log \diam D_n = O( |D_n|^\delta),
\forall \delta > 0$, then
$$|D_n|^{-\frac12} \, \sum_{\el \in D_n} \, T^\el f(.)
\overset{distr} {\underset{n \to \infty} \longrightarrow } \Cal
N(0,\sigma^2(f)).$$
\end{thm}
\proof The sequence $w_n(\el) = 1_{D_n}(\el)$ is $\zeta$-regular,
with $\zeta = \delta_0$. Suppose that $\varphi_f(0) \not = 0$. We
have $\sigma_n^2(f) \sim |D_n| \, \varphi_f(0)$ and $w_n(\el) = 0$ or 1. Condition (\ref{cumulCNab}) reads here
$$\sum_{(\el_1, ..., \el_r) \in {\cal N}(\Phi,r)} \prod_{j= 1}^r 1_{\el_j \in D_n} = o(|D_n|^{r\over 2}), \text{ for } r \geq 3.$$
For $r \geq 3$, by Proposition \ref{bndzerocuml} we have
\begin{eqnarray*}
\# \{\uA \in D_n^r : C_{\widetilde P}(\uA) \not = 0\} = O(|D_n|^{{r\over 2} - \frac12} \, (\log \diam D_n)^{\theta(r)}).
\end{eqnarray*}
By the hypothesis on the diameter, this bound implies (\ref{cumulCNab}) and the result follows from Theorem \ref{tclRegKer}. \eop

{\it Remark}: For the case of rectangles, see Theorem \ref{functLim}.

\vskip 3mm
{\bf Example 2: Random walks and quenched CLT} \label{applRW}

Using the notations and results of \cite{CohCo15}, now we apply the previous sections to random walks of commuting
endomorphisms or automorphisms on a shift-invariant subgroup $G$.

Let us present the result for $d=2$. We take two polynomials $(R_1, R_2)$ with $\gamma_{R_1}, \gamma_{R_2}$ the corresponding endomorphisms
of $K = \F_p^{\Z}$ generating a 2-dimensional action with Lebesgue spectrum.
Taking the natural invertible extension, we extend them to automorphisms (the shifts $\sigma_{1}, \sigma_{2}$)
of the shift-invariant subgroup $G_{\Cal J }$ of $G_0^{(3)}$ defined by the ideal $\Cal J$ generated in $\cal P_{3}$
by $x_{2} -R_1(x_1), x_{3} -R_2(x_1)$.

Let $(X_k)_{k \in \Z}$ be a sequence of i.i.d. $\Z^2$-valued random
variables generting a reduced aperiodic random walk.

\begin{thm} \label{main} Suppose that $W$ has a finite moment of order 2 on $\Z^2$.
Let $\el \to T^\el$ be a $\Z^2$-action generated by shifts $\sigma_2, \sigma_3$  on $G_{\cal J}$. Let $f$ be in $AC_0(G_{\cal J})$
with spectral density $\varphi_f$ such that $\varphi_f(0) \not = 0$. Then, there exists a constant $C$ such that, for a.e. $\omega$,
$$\displaystyle (C n \Log n)^{-\frac12} \sum_{k=0}^{n-1} \, T^{Z_k(\omega)} f(.) \ \overset{distr} {\underset{n \to \infty} \longrightarrow }
\ \Cal N(0, 1).$$
\end{thm} \proof Theorem 4.16 in \cite{CohCo15} gives the $\delta(0)$-regularity for the r.w. summation
$(w_n(\omega, \el))_{n \geq 1} = (\sum_{k=0}^{n-1} 1_{Z_k(\omega) = \el})_{n \geq 1}$.

For a recurrent 2-dimensional r.w., for a.e. $\omega$, $\sum_\el w_n^2(\omega, \el) \sim \E \sum_\el w_n^2(., \el) \sim C n \Log \, n$.
To concluded, we need the bound:
\begin{eqnarray}
\sum_{(\el_1, ..., \el_r) \in {\cal N}(\Phi,r)} \prod_{j= 1}^r w_n(\el_j) = o((n \Log \, n)^{r/2}). \label{majcumul1}
\end{eqnarray}
For every $\delta >0$, by the law of iterated logarithm there is a finite constant $C(\omega)$ such that
$$\|\el\| > C(\omega) n^{\frac12 + \delta} \Rightarrow w_n(\omega, \el) = 0.$$
Therefore, the previous sum can be restricted to $\el_j$ in a ball of radius $C(\omega) n^{\frac12 + \delta}$.
Moreover, we know that $\sup_\el w_n(\el) = o(n^\varepsilon), \forall \varepsilon > 0$ (Proposition 4.1, in \cite{CohCo15}).
It follows that the lhs of (\ref{majcumul1}) is less than $n^{r \varepsilon}$ multiplied by the cardinal of $r$-tuples in the set
${\cal N}(\Phi,r)$ supported in the ball $B(0, C(\omega) n^{\frac12 + \delta})$, for which a bound is given by Proposition \ref{bndzerocuml}.

This bound is less than $C \, n^{r \varepsilon} \,  n^{2 \, (\frac12 + \delta) \, (\frac r2 - \frac12)}$, up to a logarithmic factor.
Taking into account only the powers of $n$, on the left hand, we find for the power of $n$:
${r\varepsilon} + (1 + 2\delta) \, (\frac r2 - \frac12)$ which is $< r/2$, if $\varepsilon +\delta < \frac12 r^{-1}$.
\eop

\section {\bf Appendix: endomorphisms of compact abelian groups} \label{append0}

We recall here some properties of endomorphisms of compact abelian groups.
$G$ denotes a compact abelian group, $\widehat G$ the dual group of characters on $G$, $\chi_{0}$ the trivial character.

The following fact has been used in the first section: let $H$ be a closed subgroup of $G$, $L$ a subgroup of $\widehat G$.
If $H^\perp = \{\chi \in  \widehat G: \chi(h) = 1, \, \forall h \in H \}$ denotes the subgroup
of $\widehat G$ annulator of $H$ and if $L^\perp = \{h \in  G: \chi(h) = 1, \, \forall \chi \in L \}$ denotes the closed subgroup
of $G$ annulator of $L$, then $(H^\perp)^\perp = H, \ (L^\perp)^\perp = L$.

Let $d \geq 1$ be an integer and $(T_1, ..., T_d)$ commuting endomorphisms of $G$.
If $\el=(\ell_1, ..., \ell_d)$ is in $\N^d$, we write $T^\el$ for $T_1^{\ell_1}... T_d^{\ell_d}$.
If $f$ is function on $G$, $T^\el f$ stands for $f \circ T^\el$.

If necessary, we lift the action to an invertible $\Z^d$-action by commuting automorphisms of an extension of $G$.
The $\Z^d$-action is said to be {\it totally ergodic} if $T^\el$ is ergodic for every $\el \in \Z^d \stm0$.
It is equivalent to: $T^\el \chi \not = \chi$ for  $\el \not = \underline 0$ and any character $\chi \not = \chi_0$,
to the Lebesgue spectrum property, as well as to $2$-mixing.

Let $(f_1, ..., f_r)$ be a finite set of trigonometric polynomials and $\Phi = (\chi_j, j \in J)$ be the finite set of characters $\not = \chi_0$
on $G$ such that $f_i = \sum_{j \in J} c_{i,j}(f_i) \, \chi_j$, $j=1, ..., r$.
For $(\a_1, ... \a_r) \in (\Z^d)^r$, we have
\begin{eqnarray}
&&\int f_1(T^{\a_1} x) ... f_r(T^{\a_r} x) \ dx =
\sum_{j_1, ..., j_r \in J} c_{1, j_1} ... c_{r, j_r} 1_{T^{\a_1}
\chi_{j_1} ... T^{\a_r} \chi_{j_r} = \, \chi_{0}}. \label{r-moment}
\end{eqnarray}

{\it Exactness}

If $\gamma$ is a surjective algebraic endomorphism of $G$, its action on $\hat G$, still denoted by $\gamma$, is injective.
The operator of composition by $\gamma$ on $L^1(G)$ is denoted by $T_\gamma$.
In what follows we consider endomorphisms with {\it finite kernel}.

The adjoint operator $\Pi_\gamma$ of $T_\gamma$ is defined by
$$\int_G T_\gamma f \ g \, d\mu = \int_G \, f \ \Pi_\gamma g \, d\mu, \ f \in L^1, g \in L^\infty.$$
It is a contraction of $L^\infty(G)$ and it extends to a contraction of $L^2(G)$.
It can be expressed for $f = \sum_{\chi \in \hat G} \ c_f(\chi) \, \chi(x) \in L^2(G)$ as
\begin{eqnarray}
\Pi_\gamma f (x) &=& {1 \over |K_\gamma|} \sum_{y: \, \gamma y = x} f(y) = \sum_{\chi \in \hat G} \ c_f(\gamma \chi) \, \chi(x)  \label{defP2}.
\end{eqnarray}
It follows from (\ref{defP2}) that
\begin{eqnarray}
\Pi_\gamma \, \chi_1 = 0 \text{ if } \chi_1 \not \in \gamma \hat G, \ \Pi_\gamma \, \chi_1 = \chi_2
\text{ if there is } \chi_2 \in \hat G \text{ such that } \gamma \chi_2 = \chi_1. \label{Pcharac}
\end{eqnarray}
By injectivity, $\chi_2$ is uniquely defined in the second case.

Recall that an endomorphism $\gamma$ is exact (as a measure preserving map on $(G, \mu)$), if
\begin{eqnarray}
\lim_n \|\Pi_{\gamma}^n f\|_2 = 0, \ \forall f \in L_0^2(\mu). \label{exact1}
\end{eqnarray}
Exactness of an (algebraic) endomorphism $\gamma$ is equivalent to:
\begin{eqnarray}
\forall \chi \not = \chi_0, \ \exists N(\chi) \text{ such that }
\Pi_\gamma^n \chi = 0, \text{ for } n \geq N(\chi). \label{exact2}
\end{eqnarray}

Let $R$ be in $\F_p[x]$. It defines an endomorphisms $\gamma_R$ of the group $\F_p^{\Z^+}$.
By (\ref{defP2}), the transfer operator $\Pi = \Pi_{\gamma_R}$ acts on a function formally defined by its Fourier series as follows:
$$f = \sum_{P \in \F_p[x]} c(f, \chi_{P}) \, \chi_P  \to \Pi f = \sum_{P \in \F_p[x]} c(f, \chi_{PQ}) \chi_Q.$$
Therefore we have: $\Pi^n f = \sum_Q c(f, \chi_{R^nQ}) \chi_Q$ and $\|\Pi^n f\|^2 = \sum_Q |c(f, \chi_{R^nQ})|^2$.

A character $\chi= \chi_P$ associated to a polynomial $P \in \F_p[x^\pm]$,
belongs to $R(x)^n \F_p[x^\pm]$ if $P$ is divisible by $R(x)^n$. Therefore, either $R(x) = c x^\varepsilon$,
with $c \not = 0$ in $\F_p$ and $\varepsilon = \pm 1$ or $\gamma_R$ is {\it exact},
since then, for every a polynomial $P$, there is $N(P)$ such that $P$ is not divisible by $R(x)^n$ for $n \geq N(P)$.

{\it Complete commutation}
\begin{prop} \label{totCommEquiv} Let $\gamma_1, \gamma_2$ be commuting surjective endomorphisms of $G$ such that $\Ker(\gamma_1)$ is finite.
The following conditions are equivalent \,\footnote{\ Property (\ref{totComm}) is the notion of complete commutation used by M.~Gordin \cite{Go09}.
See also \cite{CuVe12}. This property can be viewed as a primality condition between $\gamma_1$ and $\gamma_2$.}:
\begin{eqnarray}
T_{\gamma_2} \Pi_{\gamma_1} &=&  \Pi_{\gamma_1}  T_{\gamma_2}, \label{totComm} \\
\gamma_2 \chi \in \gamma_1 \hat G &\Rightarrow& \chi \in \gamma_1 \hat G, \label{gamma12}\\
\Ker(\gamma_1) &\cap& \Ker(\gamma_2) = \{0\}.  \label{Ker12}
\end{eqnarray}
\end{prop}
\proof
Condition (\ref{totComm}) is equivalent to $T_{\gamma_2} \Pi_{\gamma_1} \chi =  \Pi_{\gamma_1}  T_{\gamma_2} \chi$, for every $\chi \in \hat G$.

Using (\ref{defP2}), we have $T_{\gamma_2} \Pi_{\gamma_1} \chi = 0$ if $\chi \not \in \gamma_1 \hat G$, $= \gamma_2 \zeta$
if $\chi = \gamma_1 \zeta$, with $\zeta \in \hat G$.
Likewise, we have $\Pi_{\gamma_1} \gamma_2 \chi = 0$ if $\gamma_2 \chi \not \in \gamma_1 \hat G$, $= \eta$
if $\gamma_2 \chi = \gamma_1 \eta$, with $\eta \in \hat G$.

Therefore, (\ref{totComm}) is equivalent to: $\gamma_2 \chi \not \in \gamma_1 \hat G \Leftrightarrow \chi \not \in \gamma_1 \hat G$,
i.e., to (\ref{Ker12}), since the implication $\Leftarrow $ is always satisfied by commutativity.

The annulator of $\gamma_1 \hat G$ is the kernel of $\gamma_1$. By commutation of $\gamma_1$ and $\gamma_2$, the kernel $Ker(\gamma_1)$ is
mapped into itself by $\gamma_2$.

By (\ref{gamma12}), $\Ker(\gamma_1)$ and $\gamma_2 \Ker(\gamma_1)$ have the same annulator, hence they coincide.
The equality $\gamma_2 \Ker(\gamma_1) = \Ker(\gamma_1)$ implies that $\gamma_2$ is surjective on $\Ker \gamma_1$.
Since the kernel is finite, injectivity and surjectivity of the restriction of $\gamma_2$ to $K(\gamma_1)$
are equivalent. Therefore, injectivity holds.

Now let $u \in \Ker(\gamma_1) \cap \Ker(\gamma_2)$. It satisfies $u \in \Ker(\gamma_1)$ and $\gamma_2 u = 0$.
By injectivity of the restriction of $\gamma_2$ to $Ker(\gamma_1)$, this implies $u = 0$.

Conversely, the condition $\Ker(\gamma_1) \cap \Ker(\gamma_2) = \{0\}$ implies injectivity, hence surjectivity and (\ref{gamma12}) follows.
\eop

The symmetry in Condition (\ref{Ker12}) implies: $T_{\gamma_1} \Pi_{\gamma_2} =  \Pi_{\gamma_2}  T_{\gamma_1}$.
Observe that if $\gamma_1 = \gamma_2$ the equivalent conditions are satisfied if and only if $\gamma_1$ is an automorphisms.

{\it Example 1}: (Endomorphisms of $\T^\rho$, $\rho > 1$) Let $A, B$ be two commuting non singular matrices $d \times d$ with coefficients in $\Z$.
A sufficient condition for (\ref{Ker12}) for the endomorphisms defined by $A$ and $B$ on $\T^\rho$ is that, in the decomposition of $\R^d$
into irreducible (over $\Z$) spaces $V_j$ under $A$ (and $B$), for each $V_j$ the determinants of the restriction of $A$ and $B$ are relatively prime.
See also \cite{CuDeVo16}.

{\it Example 2}: (Endomorphisms of $\F_p^{\Z^+}$)\  If $R_1, R_2$ are two relatively prime polynomials in one variable,
the endomorphisms $\gamma_{R_1}$ and $\gamma_{R_2}$ acting on $\F_p^{\Z^+}$ endowed with its Haar measure are completely commuting.
This follows from Bezout relation and (\ref{gamma12}).

\vskip 3mm
{\it Regular functions on $\F_p^{\Z^d}$}

A distance $\rho$ on $\F_p^{\Z^d}$ is defined by $\rho(\zeta, \zeta') = \sum_{\k \in \Z^d} \, 2^{-\|\k\|} |\zeta_\k - \zeta'_\k|$.
Let $D_n$ be the square  $\{\el: |\ell_1| \leq n, ..., |\ell_d| \leq n\}$. The regularity of a function $f$ on $G_0$ or on a subset of $G_0$
(how it depends on the remote coordinates) is measured by the variations
$$V_n(f) := \sup_{\zeta,\zeta': \, \zeta_{\el} = \zeta'_{\el}, \, \forall \el \, \in D_n} \, |f(\zeta) - f(\zeta')|, \ n \geq 1.$$
If $D$ is a finite set in $\Z^d$, the space $\Cal F(D)$ of complex valued functions on the finite group $\F_p^{D}$
can be viewed as the subspace of the space $\Cal C(\F_p^{\Z^d})$ of complex continuous functions on $\F_p^{\Z^d}$ depending
only on the coordinates $\zeta_\el$ for $\el \in D$.

To a point $\zeta$ in $\F_p^{\Z^d}$, let us associate the point $\pi_D(\zeta)$ whose coordinates coincide with the coordinates of $\zeta$ on $D$
and are equal to 0 outside $D$. If $f$ is a function on $\F_p^{\Z^d}$, we denote by $\Pi_D f$ the function in $\Cal F(D)$ defined by
$\Pi_D f(\zeta) = f(\pi_D(\zeta))$.

If $(D_n)_{n \geq 0}$ is an increasing sequence of domains in $\Z^d$ such that $\bigcup_n D_n = \Z^d$, then for every continuous function $f$
on $\F_p^{\Z^d}$ we have: $\lim_n \|f - \Pi_{D_n} f \|_\infty = 0$. It follows that $f = \sum_{n= 0}^\infty \, \varphi_n$,
where $\varphi_n = \Pi_{D_n} f - \Pi_{D_{n-1}} f\in \Cal F(D_n)$ and the series is converging in the uniform norm.

An approximation of $f$ depending only on coordinates in $\{0, ..., n-1\}$ is yield by replacing $f$ by the function $\varphi_n$
such that $\varphi_n(\zeta)= f(\pi_n(\zeta))$. Clearly we have $\|f - \varphi_n\|_\infty \leq V_n(f)$.

Suppose that $\varphi_n$ is a function on $\F_p^{\Z_+}$ depending only on coordinates in $\{0, ..., n-1\}$. Then $\varphi_n$ is a finite sum
of characters supported on subsets of $\{0, ..., n-1\}$. If the polynomial $R$ has degree at least 1, $R^nQ$  has degree $\geq n + \deg(Q)$.
Polynomials of degree respectively $\geq n + \deg(Q)$ and $< n$ define orthogonal characters. It follows that $\varphi_n$
is orthogonal to characters of the form $\chi_{R^n Q}$.

\begin{lem} \label{regul1} Let $f$ satisfy $V_n(f) =O(\lambda^n)$, for $\lambda < p^{-1}$. Then $f$ belongs to $AC_0(\F_p^{\Z_+})$
and for every $R \in \Cal P[x]$ of degree $\geq 1$, $\|\Pi^n f\|_\infty \leq C'\lambda^{n}$, for a constant $C$.
\end{lem}
\proof We have:
$\sum_{R}|c(\chi_Q, f)| \leq \sum_n \sum_{deg(Q) < n} \, |c(\chi_Q, f)| \leq \sum_{n \geq 1} \#\{Q: \deg Q < n\} \, V_{n}(f)
= O(\sum_{n \geq 1} p^n \lambda^n)= O(1)$.

Writing $f = f - \varphi_{n+k} + \varphi_{n+k}$, we obtain
$c(f, \chi_{R^nQ}) = \langle f, \chi_{R^n Q}\rangle = \langle f - \varphi_{n+k}, \chi_{R^n Q} \rangle$, if $k < \deg Q$,
which implies $|c(f, \chi_{R^nQ})| \leq \|f - \varphi_{n+k}\|_\infty \leq V_{n+k}(f)$.

We deduce that $\|\Pi^n f\|_\infty \leq \sum_Q |c(f, \chi_{R^nQ})| \leq \sum_Q V_{n+\deg Q}(f)$ is bounded by
\begin{eqnarray*}
\sum_{k \geq 0} \#\{Q: \deg Q = k\} \, V_{n+k}(f) \leq \sum_{k \geq 0} \, p^{k+1} \, V_{n+k}(f) \leq C\lambda^{n}. \eop
\end{eqnarray*}

{\it A limit theorem for sums of rectangles}

The case of rectangles is a special case
for which the martingale method can be used to obtain a functional theorem for ergodic sums.
From Theorems 1 and 8 in \cite{CuDeVo16} and the previous lemma, it can be deduced:
\begin{thm} \label{functLim} Let $(R_j, j= 1,..., d)$ be pairwise relatively prime polynomial of degree $\geq 1$ and let
$(\gamma_j= \gamma_{R_j})$ be the associated family of commuting algebraic exact endomorphisms of $\F^{\Z^+}$,
such that $\Ker(\gamma_i) \cap \Ker(\gamma_j) = \{0\}$, for $i \not = j$. Let $T_j := T_{\gamma_j}, j= 1,..., d$
and $T^\el = T_1^{\el_1}... T_d^{\el_d}$. If $(D_n)_{n \geq 1}$ is an increasing sequence of rectangles and if
$f$ satisfies the regularity condition $V_n(f) = O(\lambda^n)$ with $\lambda < p^{-1}$, then the sequence
$(|D_n|^{-\frac12} \sum_{\el \in D_n} T^\el \, f)_{n \geq 1}$ satisfies a functional CLT.
\end{thm}

\vskip 3mm {\bf Acknowledgements} This research was carried out
during visits of the first author to the IRMAR at the University of Rennes 1 and
of the second author to the Center for Advanced Studies in
Mathematics at Ben Gurion University.
The first author was partially supported by the ISF grant 1/12.
The authors are grateful to their hosts for their support.

\end{document}